\documentclass[preprint,11pt]{elsarticle}
\usepackage[T1]{fontenc}
\usepackage{lmodern,amsmath,amssymb,amsthm,mathtools,booktabs,array,tikz,microtype}
\usepackage{etoolbox,needspace,enumitem}
\usepackage[hidelinks]{hyperref}
\usetikzlibrary{arrows.meta}
\newtheorem{theorem}{Theorem}[section]
\newtheorem{lemma}[theorem]{Lemma}
\newtheorem{proposition}[theorem]{Proposition}
\newtheorem{corollary}[theorem]{Corollary}
\theoremstyle{definition}
\newtheorem{definition}[theorem]{Definition}
\theoremstyle{remark}
\newtheorem{remark}[theorem]{Remark}
\newtheorem{example}[theorem]{Example}
\DeclareMathOperator{\QS}{QS}
\DeclareMathOperator{\BiLip}{BiLip}
\DeclareMathOperator{\Iso}{Iso}

\DeclareMathOperator{\id}{id}
\DeclareMathOperator{\intr}{int}
\DeclareMathOperator{\diam}{diam}
\DeclareMathOperator{\dist}{dist}
\DeclareMathOperator{\conv}{conv}

\let\C\relax
\newcommand{\C}{\mathbb C}

\newcommand{\N}{\mathbb N}
\newcommand{\Z}{\mathbb Z}
\hypersetup{pdftitle={Quasisymmetry of self-similar carpets},
pdfauthor={Xiaohua Wang; Fan Wen}}

\makeatletter
\patchcmd{\pprintMaketitle}{\hrule\vskip12pt}{\vskip12pt}{}{}
\def\ps@pprintTitle{\let\@oddhead\@empty\let\@evenhead\@empty
\let\@oddfoot\@empty\let\@evenfoot\@empty}
\makeatother
\begin{document}
\begin{frontmatter}
\title{Quasisymmetry of self-similar carpets}
\nonumnote{Supported by the National Natural Science Foundation of China (NSFC), Grant Nos. 12571094 and 1250010257}
\author[cufe]{Xiaohua Wang}
\ead{fridman\_wxh@163.com}
\author[jnu]{Fan Wen\corref{cor}}
\ead{wenfan@jnu.edu.cn}
\cortext[cor]{Corresponding author}
\address[cufe]{China Economics and Management Academy, Central University of Finance and Economics, Beijing, China}
\address[jnu]{Department of Mathematics, Jinan University, Guangzhou, China}
\begin{abstract}
We classify self-similar carpets obtained by repeating a fixed pattern
in the same orientation on a square grid. All boundary cells are
retained. The first approximation has a nonempty family of Jordan
polyomino holes with pairwise disjoint closures in the interior of the
unit square. Every quasisymmetry between two such carpets is the
restriction of a square symmetry. Thus quasisymmetric, bi-Lipschitz,
and isometric equivalence coincide. Two patterns give equivalent
carpets exactly when their subdivision parameters have a common positive
power and their retained cells agree up to symmetry on the common grid.
The criterion extends to patterns periodic from the first generation.
We also construct an area-zero square carpet with uniformly relatively
separated peripheral circles and a quasisymmetric involution that sends
a corner into an open side. Its quasisymmetry group is dihedral of
order eight, and its isometries are the four rotations. It belongs to
the quasisymmetric class of the standard middle-third carpet and to a
different isometry class.

\end{abstract}
\begin{keyword}
Sierpi\'nski carpet \sep quasisymmetric rigidity \sep self-similar set
\sep Schottky map \sep uniformization \sep isometric classification
\MSC[2020] 30C65, 28A80
\end{keyword}
\end{frontmatter}
\clearpage

\section{Introduction}

Suppose that a planar carpet lies in the unit square and contains its boundary.
Its four corners are exactly the endpoints of pairs at distance
$\sqrt{2}$. An isometry between two such carpets must therefore
send corners to corners. A quasisymmetry controls ratios of distances
\cite{TV80}. Under what conditions does it still distinguish a corner
from an interior point of a side?

We study carpets for which each retained cell contains a scaled
copy of the same deletion pattern in the same orientation. Every
boundary cell is retained. Adjacent deleted cells may form one hole,
so that an actual hole can be a nonconvex polygon. We require these
holes to have Jordan boundaries and pairwise disjoint closures in the
interior of the square. In this class the grid forces a quasisymmetry
to recognize the corners and determines its values on the entire
carpet. The resulting classification compares different subdivision
parameters as well as different deletion patterns.

Example~\ref{scope:markedexample} gives a different answer to the
corner question for a carpet with uniformly relatively separated
square holes in the same unit square.

\subsection{The class of carpets and the main theorem}

Let $Q=[0,1]^2$, $U=\intr Q$, and $O=\partial Q$.
Fix an integer $n\ge3$ and a proper digit set
$D\subsetneq\{0,\ldots,n-1\}^2$ satisfying
\begin{equation}\label{main:frame}
 \{(i,j):i\in\{0,n-1\}\ \text{or}\ j\in\{0,n-1\}\}\subset D.
\end{equation}
For $d\in D$, put $\varphi_d(z)=(z+d)/n$, and let $K=K(n,D)$ be the
unique nonempty compact set such that
\begin{equation}\label{main:ifs}
 K(n,D)=\bigcup_{d\in D}\varphi_d(K(n,D)).
\end{equation}
All copies have the same orientation. We call~\eqref{main:frame} the
full frame condition: every boundary cell is retained.

We use the closed retained-cell convention. The first approximation is
$Q_1=\bigcup_{d\in D}\varphi_d(Q)$, and its actual holes are the
components of $U\setminus Q_1$. A common side of two deleted cells is
not added back to the set. We assume
\begin{equation}\tag{P}\label{main:P}
 \begin{gathered}
 U\setminus Q_1=\bigcup_{j=1}^{h}V_j,\qquad h\ge1,\\
 V_j\text{ is a Jordan polyomino domain},\qquad \overline{V_j}\subset U,\\
 \overline{V_i}\cap\overline{V_j}=\varnothing\quad(i\ne j).
 \end{gathered}
\end{equation}
Here a Jordan polyomino domain is an open Jordan domain obtained by
merging finitely many grid cells; its boundary is an axis-parallel simple
polygon. Condition~\eqref{main:P} allows adjacent deleted cells, but
excludes contact between distinct actual holes and self-contact of a
hole boundary. All carpets carry the restricted Euclidean metric.

Recall that a homeomorphism $f:X\to Y$ is quasisymmetric if there is an
increasing homeomorphism $\eta:[0,\infty)\to[0,\infty)$ such that
\[
 \frac{|f(x)-f(y)|}{|f(x)-f(z)|}
 \le\eta\!\left(\frac{|x-y|}{|x-z|}\right)
\]
for all distinct $x,y,z\in X$. We write $\QS(X,Y)$ for the set of such
homeomorphisms onto $Y$, and $D_4=\Iso(Q)$ for the group of the eight
Euclidean symmetries of $Q$.

\begin{theorem}\label{main:classification}
Let $K_1=K(n,D)$ and $K_2=K(m,E)$, where $n,m\ge3$, and suppose that
both constructions satisfy~\eqref{main:frame} and~\eqref{main:P}. Then
\begin{equation}\label{main:allmaps}
 \QS(K_1,K_2)
 =\{R|_{K_1}:R\in D_4,\ R(K_1)=K_2\}.
\end{equation}
If this set is nonempty, there are positive integers $a,b$ such that
$n^a=m^b$. For any common power $N=n^a=m^b$, define
\begin{align}
 D^{[a]}&=\left\{\sum_{j=1}^{a}n^{a-j}d_j:d_j\in D\right\},
 \label{main:digitsD}\\
 E^{[b]}&=\left\{\sum_{j=1}^{b}m^{b-j}e_j:e_j\in E\right\}.
 \label{main:digitsE}
\end{align}
For $R(z)=A_Rz+c_R\in D_4$, let
\begin{equation}\label{main:rho}
 \rho_N^R(d)=A_Rd+(N-1)c_R.
\end{equation}
Then
\begin{equation}\label{main:finite}
 R(K_1)=K_2\quad\Longleftrightarrow\quad
 E^{[b]}=\rho_N^R(D^{[a]}).
\end{equation}
Thus the two carpets are quasisymmetrically equivalent if and only if
their subdivision parameters have a common power and their grouped
digit sets differ by a square symmetry.
\end{theorem}

\begin{corollary}\label{main:isometricclassification}
Under the hypotheses of Theorem~\ref{main:classification},
\begin{equation}\label{main:threeclasses}
 \QS(K_1,K_2)=\BiLip(K_1,K_2)=\Iso(K_1,K_2),
\end{equation}
where $\BiLip$ denotes the bi-Lipschitz homeomorphisms and $\Iso$ the
surjective isometries. For a single carpet in this class,
\begin{equation}\label{unified:selfgroup}
 \QS(K)=\BiLip(K)=\Iso(K)
 =\{R|_K:R\in D_4,\ R(K)=K\}.
\end{equation}
For fixed $n$, the quasisymmetric and isometric classes are both
parametrized by $\mathfrak D_n/D_4$, where $\mathfrak D_n$ is the finite
set of digit sets satisfying~\eqref{main:frame} and~\eqref{main:P},
with the action $D\mapsto\rho_n^R(D)$.
\end{corollary}

\begin{proof}
Every isometry is bi-Lipschitz. If $f$ is $L$-bi-Lipschitz, then
\[
 \frac{|f(x)-f(y)|}{|f(x)-f(z)|}
 \le L^2\frac{|x-y|}{|x-z|},
\]
so it is quasisymmetric with control function $\eta(t)=L^2t$.
Hence $\Iso\subset\BiLip\subset\QS$, and the reverse inclusion
$\QS\subset\Iso$ follows from Theorem~\ref{main:classification}.
Taking $K_1=K_2$ gives~\eqref{unified:selfgroup}. For fixed $n$,
take $a=b=1$ in~\eqref{main:finite} to obtain the orbit description.
\end{proof}

Fixing $Q$ removes the freedom to rescale. Grouping two ternary
generations still gives a $9\times9$ description of the same set.
For separated unit-cell holes the largest
hole recovers the original parameter (Section~\ref{unified:applications}).

In Example~\ref{scope:markedexample}, four points in the interiors of
the sides of the standard middle-third carpet $S_3$ become the corners
of a new square carpet $T$. Under this change of parametrization the
quarter-turn remains a rotation, while a reflection moves a new corner
into an open side. Conjugation gives all eight quasisymmetries of $T$.
Only the four rotations preserve distance. Since $S_3$ has eight
isometries, no isometry can carry $S_3$ onto $T$.

\subsection{Related work}

Bonk and Merenkov~\cite{BM2013} proved rigidity of the standard
middle-third carpet $S_3$. For the central-hole carpets $S_p$, with
$p\ge3$ odd, they distinguished the quasisymmetric classes using
carpet modulus and weak tangents. Their carpet modulus is related to
Schramm's transboundary extremal length~\cite{Sch95}. Rao and
Shengyou Wen~\cite{RW18} proved rigidity of $S_5$, and Bonk and
Merenkov~\cite{BM2020} proved it for every $S_p$. Zeng and
Su~\cite{ZS15} classified a family of symmetric four-hole carpets.
Rao and Fan Wen~\cite{RW25} established finiteness of the
quasisymmetry group for balanced self-similar quasi-round carpets.

For general square-grid digit sets, Lau, Luo and Rao~\cite{LLR13}
used the periodic extension $K+\Z^2$ to classify connected-component
structures. Ruan and Wang~\cite{RuanW17} classified ternary fractal
squares with one or two deleted cells up to bi-Lipschitz equivalence.
Dai, Luo, Ruan, Wang and Xiao~\cite{DLRWX23} gave criteria for
connectedness and for the occurrence of cut points and local cut points.
For a connected full-frame construction with a nonempty deletion,
Proposition~\ref{geo:nolocalcut} identifies the absence of local cut
points with condition~\eqref{main:P} on the first-generation holes.

Our proof uses Bonk's round uniformization~\cite{Bonk2011} and the
local theory of Schottky maps. The latter builds on the M\"obius
rigidity of measure-zero Schottky sets~\cite{BKM09} and on
Merenkov's regularity, uniqueness and stabilization
results~\cite{Mer12,Mer14}. We use the local statements in
\cite[Section~4]{BM2020}. Ntalampekos~\cite{Nta26} characterized
quasiconformal equivalence to Schottky sets under a condition on pairs
of complementary components. Related continuation arguments occur in
work on carpet Julia sets by Bonk, Lyubich and Merenkov~\cite{BLM16}
and by Merenkov and Shen~\cite{MS26}. Merenkov~\cite{Mer26} studied
local distinctions between Kleinian carpets and rational carpet
Julia sets.

The construction of $T$ uses Ntalampekos's square
uniformization~\cite[Chap.~3]{Nta20}, in which four boundary arcs
determine the sides of the target rectangle. The four-corner theorem
of Bonk and Merenkov~\cite[Thm.~1.4]{BM2013} determines its aspect
ratio. They relate this uniqueness statement to Schramm's theorem for
rectangular square tilings with prescribed combinatorics~\cite{Sch93}.
For an area-zero square carpet in a rectangular ring, Fan
Wen~\cite{Wen25} proved that a quasisymmetric self-map fixing all
four outer vertices is the identity.

In the broader setting of metric carpets, Cheeger and
Eriksson-Bique~\cite{CEB23} constructed infinitely many pairwise
quasisymmetrically inequivalent Loewner carpets with square-carpet
realizations. Sheng, Wei and Wen~\cite{SWW24} constructed rigid
metric carpets of arbitrarily large Hausdorff dimension with
rectifiable peripheral circles. Merenkov~\cite{Mer10} gave a carpet
for which every quasisymmetric self-embedding is onto. On the carpets
of Theorem~\ref{main:classification}, each retained-cell similarity
is a quasisymmetric embedding onto a proper subset.

\subsection{The mechanism of the proof}

Normalizing on selected cells gives local maps $H_i$ with a
nondegenerate quasiconformal limit $H$. Convergence of the maps and
their inverses preserves the exact image identities. Schottky
stabilization applied to $H_iH^{-1}$ on the target carpet gives
agreement on a nonempty relatively open set.

The propagation step concerns a conjugated cell contraction that agrees
with a target-grid homothety on this set. A common open side of two
retained cells has no peripheral points. Its image therefore forces
the adjacent target cell to be retained. Only after this is established
do we use local exactness of the homothety and Schottky uniqueness to
continue the identity. The result is a relation
\[
 fA^k=Lf\quad\text{on }K_1,
\]
where $A$ is a source contraction and $L$ is a positive target
homothety. The argument allows one exceptional peripheral circle and
does not assume preservation of the outer square.

The conjugacy relation first identifies the outer square and then its
four corners. In a thin strip along a side, translations between
retained cells give exact affine increment identities for the boundary
parametrization. Their scale factors belong to the discrete group
$m^{\Z}$ and tend to $1$, so the increments are eventually translations.
A difference-quotient squeeze determines the whole boundary map,
without a differentiability assumption. Extending by integer
translations gives a relation between global expanding maps.
Comparison with the integer periods forces $n^k=m^\ell$; bounded
periodic error then forces the original map to be an isometry.

Sections~2--4 establish the local results used in Sections~5--8 to
prove rigidity. Section~9 derives the finite criterion and treats
periodic patterns. Section~\ref{unified:applications} computes the
groups in the examples.

For a retained word $w=d_1\cdots d_k$, we write
$A_w=\varphi_{d_1}\circ\cdots\circ\varphi_{d_k}$ and call $A_w(Q)$
a retained cell. A power $A^r$ means repetition of one fixed
contraction, as distinct from an arbitrary retained word of length $r$.


\section{Geometry and connectivity}

The results of this section apply separately to $K_1$ and $K_2$.
To distinguish their subdivision parameters, we work with a generic
carpet $K=K(q,\mathcal F)$, where $q\geq3$ and
$\mathcal F\subsetneq\{0,\ldots,q-1\}^2$ contains every boundary
digit. Write
\[
 Q=[0,1]^2,\qquad U=\operatorname{int}Q,\qquad O=\partial Q,
 \qquad \phi_d(z)=\frac{z+d}{q}.
\]
For a finite word $w$ over $\mathcal F$, let
$A_w=\phi_{w_1}\circ\cdots\circ\phi_{w_{|w|}}$; the empty word
is allowed. The squares $A_w(Q)$ are the retained cells, and
$Q_k=\bigcup_{|w|=k}A_w(Q)$.
Condition~\textup{(P)} requires the nonempty components
$V_1,\ldots,V_h$ of $U\setminus Q_1$ to be Jordan polyomino domains
such that
\[
 \overline{V_j}\subset U,\qquad
 \overline{V_i}\cap\overline{V_j}=\varnothing\quad(i\ne j).
\]
Holes always mean complementary components. In particular, the common
seams of adjacent deleted cells are not retained as additional subsets.

\begin{lemma}[Geometry and exact copies]\label{geo:geometry}
Under condition~\textup{(P)}, the set $K$ is a planar Sierpi\'nski
carpet of area zero and has no local cut points. Its bounded
complementary components are precisely
\begin{equation}\label{geo:holes}
 A_w(V_j),\qquad w\in\mathcal F^*,\quad 1\leq j\leq h.
\end{equation}
The peripheral circles, including $O$, are uniform quasicircles and
are uniformly relatively separated. The carpet $K$ is locally porous.
At each level, the graph of retained cells joined across common sides
is connected. For every cell contraction $A=A_w$,
\begin{equation}\label{geo:native}
 K\cap A(Q)=A(K),\qquad K\cap A(U)=A(K\setminus O).
\end{equation}
Moreover, a $q^{-k}$ grid cell $P\subset Q$ is retained if and only if
$K\cap\operatorname{int}P\ne\varnothing$.
\end{lemma}

\begin{proof}
We first prove side connectivity at the first level. The open set
$U\setminus\bigcup_j\overline{V_j}$ is path connected. Indeed, a
polygonal path can be diverted around each closed polygonal disk along
a sufficiently small exterior polygonal collar. The disks form a
finite family, have positive mutual distances, and have positive
distance from $O$. The diversions can therefore be chosen to avoid the
other disks and the outer boundary. Join the centers of any two retained
cells by such a path. A small perturbation avoids all first-level grid
vertices and prevents the path from following a grid edge. Every change
of cell then occurs across a common open side. Hence the side-adjacency
graph is connected.

At the next level, the same connected graph replaces each retained
cell. Along a common side of two parent cells, the full boundary rows
of child cells are retained on both sides. Induction proves side
connectivity at every level. Thus each $Q_k$ is connected, and so is
the decreasing intersection $K=\bigcap_kQ_k$.

Every point of $Q\setminus K$ has a unique first deletion level.
Consequently, $U\setminus K$ is the union of the domains in
\eqref{geo:holes}. Their closures are pairwise disjoint. To see this,
compare the addresses of two holes. If the addresses first separate
at some digit, the holes lie strictly inside distinct cells at that
level. If one address extends the other, the earlier root hole is
disjoint from the interior of the retained child cell containing the
later hole, whose closure is strictly inside that child cell. Distinct
root holes in the same parent cell are separated by condition~\textup{(P)}.
These facts also show that no two such domains merge at a later stage;
they are therefore exactly the complementary components. Their closures
are disjoint from $O$.

Only finitely many holes occur up to any fixed generation, and a hole
with address $w$ has diameter at most $\sqrt2q^{-|w|}$. Thus the hole
diameters form a null family. Every retained cell contains a later hole.
Since retained-cell diameters tend to zero, $K$ has empty interior.
Whyburn's planar characterization~\cite{Whyburn1958} now shows that
$K$ is a Sierpi\'nski carpet and, in particular, has no local cut
points. If $N=|\mathcal F|<q^2$, the interiors of distinct cells of
the same level are disjoint, and hence
\[
 \mathcal L^2(K)\leq\mathcal L^2(Q_k)
       =\left(\frac{N}{q^2}\right)^k\longrightarrow0.
\]

For relative separation, put
\[
 d_0=\max_j\operatorname{diam}V_j,\qquad
 \epsilon_0=\min_j\operatorname{dist}(\overline{V_j},O)>0.
\]
If $h\geq2$, let
$\epsilon_1=\min_{i\ne j}\operatorname{dist}(\overline{V_i},
\overline{V_j})>0$; if $h=1$, set $\epsilon_1=\epsilon_0$.
For two distinct holes, choose the one with the smaller parent cell,
and denote that cell's side length by $s$. Unless the two parent cells
coincide, the other hole lies outside the interior of this smaller
parent cell. Their distance is then at least $\epsilon_0s$.
For holes in the same parent cell, the distance is at least
$\epsilon_1s$. The hole in the smaller parent cell has diameter at most
$d_0s$. Thus any two distinct inner peripheral circles $J,J'$ satisfy
\begin{equation}\label{geo:separation}
 \operatorname{dist}(J,J')\geq
 \frac{\min\{\epsilon_0,\epsilon_1\}}{d_0}
 \min\{\operatorname{diam}J,\operatorname{diam}J'\}.
\end{equation}
The same parent-cell estimate handles a pair containing $O$.

Each fixed Jordan polygon is a quasicircle. For nearby points on one
edge or on two adjacent edges, the shorter boundary arc has length
bounded by a constant times the chord length. Away from these cases,
compactness and the finite number of edges give a bounded-turning
constant. Taking the largest constant for
$O,\partial V_1,\ldots,\partial V_h$ gives one quasicircle constant
for all their positive homothetic copies.

Fix a disk $B(a,\rho)\subset V_1$, where $\rho>0$.
For $x\in K$ and sufficiently small $r>0$, choose a retained cell
$A(Q)$ containing $x$, with side length $s$ satisfying
$r/(\sqrt2q)<s\leq r/\sqrt2$.
Then $A(Q)\subset\overline B(x,r)$, and $A(B(a,\rho))$ is a
disk in the complement of $K$, of radius at least $\rho r/(\sqrt2q)$.
The diameter of its containing hole is also comparable to $r$.
This proves the required local porosity.

The full frame condition gives $A(O)\subset A(K)$.
Other cells of the same level meet $A(Q)$ only on its boundary and
therefore contribute no points outside $A(K)$. This proves
\eqref{geo:native}. An unretained grid cell has interior disjoint
from $Q_k$. Conversely, a retained cell contains an interior carpet
point: for example, $(1/q,1/q)\in K\cap U$, and its image under the
cell contraction lies in that cell's interior. This proves the cell
recognition assertion.
\end{proof}

\begin{lemma}[Common sides and grid vertices]\label{geo:edges}
The common open side of two side-adjacent retained cells of the same
level contains no peripheral point. At any grid vertex of a fixed
level, the retained cells incident to that vertex form a connected
graph when adjacency means sharing a side incident to the vertex.
\end{lemma}

\begin{proof}
Suppose the cells have level $k$. The closure of a hole created after
level $k$ lies strictly inside a retained $k$-cell and cannot meet the
common side. The boundary of an earlier hole consists of $k$-grid
edges. If it contains part of the common open side, the $k$-cell on
one side must lie in that hole, contrary to retention. Vertices of
earlier holes are $k$-grid vertices, so cannot lie in the common open
side. Finally, the open side is contained in $U$ and therefore avoids
$O$.

The only possible obstruction to the second assertion occurs at an
interior vertex where exactly two opposite quadrants are retained and
the other two are deleted. If the deleted quadrants belong to different
actual holes, their closures meet at the vertex. If they belong to
the same hole, that hole occupies two separated opposite quadrants
near the vertex and its boundary is not a simple polygon. Both cases
contradict Lemma~\ref{geo:geometry}. A vertex on the outer boundary
has only one or two adjacent quadrants inside $Q$, so the obstruction
does not arise there.
\end{proof}

\begin{lemma}[Removing finitely many peripheral circles]\label{geo:connectivity}
The set $K\setminus O$ is connected and dense in $K$. If $P,P'$ are
side-adjacent retained cells of the same level, then
\begin{equation}\label{geo:twocells}
 K\cap\operatorname{int}(P\cup P')\quad\text{is connected}.
\end{equation}
More generally, deleting any finite collection of peripheral circles
from $K$ leaves a connected dense set. The same conclusion holds for
a rectangular carpet obtained by filling every cell of a finite
rectangular array with a scaled copy of $K$ having the same orientation.
\end{lemma}

\begin{proof}
We first construct a connected skeleton. Let $\mathcal F_0$ consist
of the boundary digits of the $q$-grid alone, and put
$K_0=K(q,\mathcal F_0)\subset K$.
Let $T_j$ be the union of the boundaries of all level-$j$ retained
cells of $K_0$. The sets $T_j$ increase, and their union is dense in
$K_0$. The first-level network $T_1\setminus O$ is connected: its
inner square boundary is a connected loop, and every remaining grid
segment is attached to that loop. The same description applies when
the next network is inserted into a parent cell. If the parent cell
meets the global outer boundary, one removes one of its exterior
sides, or two adjacent exterior sides. The remaining network still
consists of the inner loop and segments attached to it, and meets
the existing network. Induction shows that every $T_j\setminus O$
is connected. Their increasing union is dense in $K_0\setminus O$;
it is also dense in $K_0$, since retained cells have diameters tending
to zero and contain network points off $O$.
Hence $K_0\setminus O$ is connected and dense in $K_0$.

Now set
\[
 E_j=\bigcup_{w\in\mathcal F^j}A_w(K_0).
\]
The sets $E_j$ increase and their union is dense in $K$.
For every $w$, the set $A_w(K_0)\setminus O$ contains the connected
set $A_w(K_0\setminus O)$ and is contained in its closure. It is
therefore connected. For side-adjacent cells these sets meet along
their common open side. Side connectivity implies that
$E_j\setminus O$ is connected. Their increasing union is contained
in $K\setminus O$ and dense in $K$, so $K\setminus O$ is connected
and dense. We have used the elementary fact that
$C\subset X\subset\overline C$ and connectedness of $C$ imply
connectedness of $X$.

For \eqref{geo:twocells}, let $e$ be the common open side.
The carpet portions in the two cell interiors are copies of
$K\setminus O$, each connected and dense in its full cell copy.
Adjoining $e$ to each preserves connectedness. The resulting sets
meet in $e$, and their union is exactly the set in
\eqref{geo:twocells}.

Let $C_1,\ldots,C_s$ be the peripheral circles to be removed; $O$
may be among them. Choose $j$ at least as large as the generation
of every inner circle in this list. Let $G_j$ be the union of all
interior carpet copies $A_w(K\setminus O)$ with $|w|=j$, together
with the common open side of every pair of side-adjacent retained
$j$-cells. Adjoining these sides to each interior copy preserves
connectedness. The connected side-adjacency graph therefore makes
$G_j$ connected. The union of the interior copies is dense in $K$,
so $G_j$ is dense as well. Every $C_i$ lies in the level-$j$ grid
skeleton and misses all the interior copies. By
Lemma~\ref{geo:edges}, the common open sides avoid all peripheral
circles. Hence
\[
 G_j\subset K\setminus\bigcup_{i=1}^sC_i
       \subset\overline{G_j}=K.
\]
The intermediate set is connected and dense.

For a finite rectangular array, first put a scaled copy of $K_0$
in each cell. Joining these along the rectangular grid gives the
same connected dense skeleton after removal of the outer rectangle.
Subdivide each cell to a common sufficiently deep level and repeat
the construction of $G_j$. All actual holes lie strictly inside
their parent cells, so common open sides again avoid the peripheral
circles. The preceding argument proves the rectangular version.
\end{proof}

\begin{proposition}[The local cut point criterion]\label{geo:nolocalcut}
Suppose that every boundary cell is retained and at least one interior
cell is deleted. If the limit set $K$ is connected, then
condition~\textup{(P)} holds if and only if $K$ has no local cut
points.
\end{proposition}

\begin{proof}
The implication from~\textup{(P)} follows from
Lemma~\ref{geo:geometry}. For the converse, assume only that $K$
is connected and has no local cut points.

First, $K$ is locally connected. For $x\in K$, let $N_j(x)$ be the
union of all level-$j$ copies $A_w(K)$ containing $x$.
These copies are connected and have the common point $x$, so
$N_j(x)$ is connected and has diameter at most $2\sqrt2q^{-j}$.
Each of the finitely many remaining level-$j$ copies has positive
distance from $x$. Thus $N_j(x)$ contains a relatively open
neighborhood of $x$. These arbitrarily small connected neighborhoods
give local connectedness.

At a first-level interior grid vertex, it is impossible to retain
exactly two opposite cells and delete the other two. Otherwise,
in a sufficiently small neighborhood of the vertex, $K$ lies in
two closed opposite quadrants meeting only at the vertex. The full
frame condition supplies line segments of $K$ approaching the
vertex in both quadrants. Local connectedness gives a connected
relatively open neighborhood of the vertex inside this small
neighborhood. Removing the vertex separates it into two nonempty
relatively open parts in the opposite quadrants. This would make
the vertex a local cut point.

The contact graph of the first-level closed retained cells is
connected. Otherwise, its finitely many components have positive
mutual distances, and each contains a nonempty copy of $K$, which
would disconnect $K$. A step in a contact path that joins only
diagonally opposite cells at a vertex can be replaced by side
steps: the excluded opposite-cell configuration guarantees another
retained cell incident to that vertex. Hence the side-adjacency
graph of retained cells is connected.

Consider the finite grid graph separating the retained and deleted
regions. It lies strictly inside $U$, because the boundary rows are
retained. At a grid vertex its degree equals the number of changes
between retained and deleted states when the four quadrants are
read cyclically. Once the opposite-cell configuration is excluded,
every vertex belonging to the boundary graph has degree two.
The graph is therefore a disjoint union of simple polygonal loops:
there are neither self-contacts nor vertices shared by two loops.

Each deleted component is consequently a polygonal domain with one
outer boundary loop and possibly finitely many inner loops.
An inner loop cannot occur. On its inner side, next to the loop,
there is a retained cell whose center can be joined to a boundary
cell by a retained-cell side path. The polygonal path joining
successive cell centers lies in the interior of the retained region
and misses the boundary graph. Yet it must cross the inner loop
to reach the outer boundary, a contradiction.
Every deleted component thus has a single boundary loop and is a
Jordan polyomino domain. The loops are pairwise disjoint, so the
closures of distinct deleted components are disjoint. They lie
strictly inside $U$ by the full frame condition. This is
condition~\textup{(P)}.
\end{proof}

\section{Round uniformization and local exactness}

\begin{definition}\label{geo:exact}
Let $X,Y$ be closed subsets of the plane, let $H$ be a homeomorphism
between ambient open sets, and let $x\in X$ lie in its domain.
We say that $H$ is \emph{locally exact from $X$ to $Y$ at $x$} if
there is an open neighborhood $W$ of $x$ such that
\begin{equation}\label{geo:exacteq}
 H(W\cap X)=H(W)\cap Y.
\end{equation}
When $X=Y$, we simply say that $H$ is locally exact for $X$ at $x$.
Both inclusions in \eqref{geo:exacteq} are required; mapping
$W\cap X$ into $Y$ alone is insufficient.
\end{definition}

\begin{lemma}[Local recognition of peripheral points]\label{geo:peripheral}
A locally exact ambient homeomorphism between planar Sierpi\'nski
carpets preserves peripheral and nonperipheral points.
\end{lemma}

\begin{proof}
A carpet point lies on a peripheral circle if and only if it is
the endpoint of an ambient embedded arc whose remaining points
lie in the complement of the carpet. One implication follows
from boundary accessibility for a Jordan domain. Conversely, the
arc without its endpoint is connected, hence lies in a single
complementary component, and the endpoint belongs to its boundary.
Shorten the arc to the neighborhood in \eqref{geo:exacteq}.
Both $H$ and $H^{-1}$ preserve this characterization.
\end{proof}

The uniform quasicircle and relative separation estimates of
Lemma~\ref{geo:geometry} permit the application of Bonk's round
uniformization theorem
\cite[Theorem~1.1 and Corollary~1.2]{Bonk2011}.
There is a quasiconformal homeomorphism $\beta$ of the sphere such
that $\beta(K)$ is a round carpet. A M\"obius normalization gives
$S=\beta(K)\subset\mathbb C$. We regard $S$ as a Schottky set
in the sphere, with its outer complementary disk included in the
family of removed disks. Quasiconformal maps preserve sets of
area zero, and local porosity is invariant under quasisymmetric
maps. Thus $S$ has area zero and is locally porous.

A quasisymmetric homeomorphism between two carpets in our class
extends to a quasiconformal homeomorphism of the sphere by
\cite[Proposition~5.1]{Bonk2011}. That proposition is stated for
quasi-M\"obius embeddings; quasisymmetric maps between the compact
carpets considered here satisfy its hypotheses. An orientation-reversing
map is allowed an orientation-reversing quasiconformal extension.

We recall the analytic results in the form needed below. A
\emph{relative Schottky set} is obtained from a region in the sphere
by removing open geometric disks whose closures are contained in
the region and are pairwise disjoint. For relative Schottky sets
$S,S'\subset\mathbb C$, a local homeomorphism
$g:W\cap S\to S'$ is a \emph{Schottky map} if
\[
 g'(z)=\lim_{\substack{w\to z\,,\;w\in W\cap S\\w\ne z}}
              \frac{g(w)-g(z)}{w-z}
\]
exists, is nonzero at every point, and is continuous on its relative
domain. This derivative is taken along the set; an ambient $C^1$
extension is not assumed.

\begin{proposition}[Schottky-map results]\label{geo:analysis}
The following are the forms of
\cite[Theorems~4.1--4.3]{BM2020} that we use. Here $S\subset\mathbb C$
is a relative Schottky set and $W$ is an ambient open set.
\begin{enumerate}
\item Suppose that $S$ has area zero and that an ambient continuous
map $g$ satisfies $g(W\cap S)\subset S$. If every point of
$W\cap S$ has a neighborhood on which $g$ is an orientation-preserving
quasisymmetric homeomorphism and is locally exact for $S$, then
$g|_{W\cap S}$ is a Schottky map. An orientation-preserving local
quasiconformal homeomorphism satisfies the quasisymmetry requirement
after the neighborhood is shrunk. Peripheral points are included.
\item Suppose that $S$ is locally porous, $W\cap S$ is connected,
and there is a Schottky map $u:W\cap S\to S$ fixing some
$a\in W\cap S$ but not equal to the identity. Let each Schottky
map $h_i$ be a homeomorphism from $W\cap S$ onto a relatively
open set $W_i\cap S$. If the maps $h_i$ converge locally uniformly
to a homeomorphism $h$ whose image is also a relatively open subset
of $S$, then $h_i=h$ on $W\cap S$ for all sufficiently large $i$.
\item Suppose that $S$ is locally porous and $W\cap S$ is connected.
If the coincidence set of two Schottky maps
$g_1,g_2:W\cap S\to S$ has a limit point in $W\cap S$, then
$g_1=g_2$ throughout $W\cap S$.
\end{enumerate}
\end{proposition}

In the stabilization statement, the nonidentity map with a fixed
point and the relatively open images are essential hypotheses.
To apply the regularity statement to the maps compared in a
uniqueness argument, local exactness must be established for both
maps; an inclusion of their images in the carpet does not suffice.
All maps to which we apply these Schottky-map results preserve
orientation. Round uniformization supplies the analytic setting,
while the grid argument below establishes the required local
exactness for successive cells.

\section{Local propagation on the target grid}

We continue to work with one fixed carpet $K=K(q,\mathcal F)$.
In applications to a map from $K_1$ to $K_2$, this section is
applied to $K=K_2$, so $q$ is the target subdivision parameter.

\begin{lemma}[Local grid propagation]\label{geo:propagation}
Let $W\subset\mathbb C$ be open, with $W\cap K$ nonempty and
connected. Let $\Phi$ be an orientation-preserving ambient
quasiconformal homeomorphism, or the restriction of an
orientation-preserving spherical quasiconformal homeomorphism
to an ambient open set avoiding its pole. Suppose that
\[
 \Phi(W\cap K)\subset K
\]
and that $\Phi$ is locally exact for $K$ at every point of
$W\cap K$. Let
\begin{equation}\label{geo:gridL}
 L(z)=q^r z+b,\qquad r\in\mathbb Z,\quad
 b\in\mathbb Z[1/q]^2.
\end{equation}
If $\Phi=L$ on a nonempty relatively open subset of $W\cap K$,
then $\Phi=L$ throughout $W\cap K$.
\end{lemma}

\begin{proof}
A Sierpi\'nski carpet is a locally path connected Peano continuum.
Thus the connected relatively open set $W\cap K$ is path connected.
The given coincidence set contains $K\cap P_0$ for some sufficiently
small retained cell $P_0$: take an ambient ball whose intersection
with $K$ lies in the coincidence set, and use the retained cells
of sufficiently small diameter. Choose a carpet point $x_0$ in
$\operatorname{int}P_0$. For an arbitrary $z\in W\cap K$, join
$x_0$ to $z$ by a path with compact image
$\gamma\subset W\cap K$. Its distance from
$\mathbb C\setminus W$ is positive.

Choose an integer $k$ sufficiently large that $k-r\geq0$,
$q^{k-r}b\in\mathbb Z^2$, and all $k$-cells meeting $\gamma$,
together with their vertex-neighbors, lie in $W$. Also require
that some retained $k$-cell $P_*$ containing $x_0$ lie entirely
in the ambient neighborhood of coincidence. Such a choice is
possible because $b\in\mathbb Z[1/q]^2$. The map $L$ then sends
the entire level-$k$ grid to the level-$(k-r)$ grid.

We replace the path by a finite cell chain without imposing
transversality to the grid. The retained $k$-cells meeting
$\gamma$ form a finite family with connected intersection graph.
Otherwise, two finite closed subfamilies would separate the
connected set $\gamma$. Hence there is a chain from $P_*$ to
a cell containing $z$, with successive cells meeting in a vertex
or a side. At a vertex-only step, Lemma~\ref{geo:edges} inserts
retained cells incident to that vertex and replaces the step by
side-adjacent steps. The inserted cells lie in $W$ by the choice
of $k$. We obtain a finite chain of side-adjacent retained cells.

Since $\Phi=L$ on $K\cap P_*$, the interior of $L(P_*)$ contains
a point of $K$. Grid alignment and the cell recognition assertion
in Lemma~\ref{geo:geometry} show that $L(P_*)$ is retained.
Suppose inductively that equality holds on $K\cap P$ for a
cell $P$ in the chain, and that $L(P)$ is retained. Let $P'$ be
the next cell and $e$ their common open side. Every point of $e$
is nonperipheral by Lemma~\ref{geo:edges}. Local exactness and
Lemma~\ref{geo:peripheral} give
\begin{equation}\label{geo:nonperipheraledge}
 \Phi(e)=L(e)\quad\text{contains only nonperipheral points of }K.
\end{equation}

We claim that $L(P')$ is retained as well. It is a grid cell of
the same level as $L(P)$ and shares a side with it. If $L(P')$
lies outside $Q$, their common open side lies in $O$, contrary
to \eqref{geo:nonperipheraledge}. If $L(P')\subset Q$ but is not
retained, its interior lies in an actual hole deleted at an earlier
level. The boundary of that hole contains the common open side
with $L(P)$, because this side belongs to the full retained
boundary of $L(P)$. Again the side is peripheral, a contradiction.
All relevant grid cells are aligned with $Q$, so there is no
case in which a cell crosses the boundary of $Q$ in its interior.
This proves the claim.

It now follows that $L$ is locally exact on the open rectangle
$B=\operatorname{int}(P\cup P')$. On both sides, the configuration
consists of two side-adjacent retained cells carrying copies of $K$
with the same orientation.
Equation~\eqref{geo:native} gives the exact copies in their
interiors, and the copies on the two sides of the common edge
fit together with the same orientation. Thus
\[
 L(B\cap K)=L(B)\cap K.
\]
The set $B\cap K$ is connected by
Lemma~\ref{geo:connectivity}. Conjugate $\Phi|_B$ and $L|_B$
by the same uniformizing map $\beta$ of the target carpet.
Part~(1) of Proposition~\ref{geo:analysis} gives two Schottky
maps. They agree on $K\cap\operatorname{int}P$ before
conjugation; this nonempty relatively open set contains a limit
point. Part~(3) of Proposition~\ref{geo:analysis} therefore
extends equality to $B\cap K$. The interior carpet is dense
in the full copy in $P'$, so continuity extends equality to
$K\cap P'$, including its boundary. Induction along the finite
chain proves $\Phi(z)=L(z)$.

In particular, no global local-exactness property of $L$ was
assumed. At every step, retention of the next target cell was
proved before the local-exactness identity was used.
\end{proof}

\begin{corollary}[Propagation with one exceptional peripheral circle]
\label{geo:onecircle}
Let $\Phi$ be an orientation-preserving quasiconformal homeomorphism
of the sphere such that $\Phi(K)\subset K$. Suppose that $\Phi$
is locally exact for $K$ at every point of $K\setminus C$, where
$C$ is a peripheral circle. If $\Phi$ agrees with a map $L$ of
the form \eqref{geo:gridL} on a nonempty relatively open subset
of $K\setminus C$, then $\Phi=L$ on $K$.
\end{corollary}

\begin{proof}
Remove from the sphere the closure of the complementary component
bounded by $C$, and also remove $\infty$ and the pole of $\Phi$.
This gives an ambient open set $W\subset\mathbb C$.
Since $\Phi(K)\subset K\subset\mathbb C$, its pole does not
belong to $K$, and hence $W\cap K=K\setminus C$.
By Lemma~\ref{geo:connectivity}, this set is connected and dense
in $K$. Apply Lemma~\ref{geo:propagation}, and then extend the
equality to $C$ by continuity.
\end{proof}


\section{Cell selection and conjugacy of contractions}\label{reset:section}

Throughout this section, let $K_1=K(n,D)$ and $K_2=K(m,E)$, and write
$Q=[0,1]^2$, $U=\intr Q$, and $O=\partial Q$. Both digit sets satisfy
the full frame condition and condition~\eqref{main:P}. Source cells have
side lengths that are negative integral powers of $n$, whereas target
cells have side lengths that are negative integral powers of $m$.
We keep these two grids separate.

Let $f\in\QS(K_1,K_2)$. By Lemma~\ref{geo:geometry} and
\cite[Proposition~5.1]{Bonk2011}, $f$ has a quasiconformal extension
\[
 F:\widehat{\mathbb C}\longrightarrow\widehat{\mathbb C},\qquad
 F|_{K_1}=f,\qquad F(K_1)=K_2.
\]
At this stage, we do not assume that $F$ preserves $O$ or $Q$, or that
$F(\infty)=\infty$.

\begin{lemma}[Cell selection in two grids]\label{reset:selection}
Let $A$ be the contraction associated with a nonempty retained word for
$K_1$. Suppose that
\begin{equation}\label{reset:inside-assumption}
 F(A^r(Q))\subset Q
\end{equation}
for all sufficiently large integers $r$. Then there are a fixed source
cell contraction $C$, a strictly increasing sequence of integers $r_i$,
and target cell contractions $B_i$ such that the maps
\begin{equation}\label{reset:Hi}
 H_i=B_i^{-1}\circ F\circ A^{r_i}\circ C
\end{equation}
have a subsequence converging uniformly in the spherical metric to a
quasiconformal homeomorphism $H$ of the same orientation. Along this
subsequence, $H_i^{-1}$ also converge uniformly in the spherical metric
to $H^{-1}$. There are constants $c_0,c_1>0$, independent of $i$, such that
\begin{equation}\label{reset:normal-bounds}
 H_i(Q)\subset U,\qquad
 \diam H_i(Q)\ge c_0,\qquad \dist(H_i(Q),O)\ge c_1,
 \qquad H(Q)\subset U.
\end{equation}
Moreover,
\begin{align}
 H_i(U\cap K_1)&=H_i(U)\cap K_2,\label{reset:exact-Hi}\\
 H(U\cap K_1)&=H(U)\cap K_2.\label{reset:exact-H}
\end{align}
Thus the normalized maps are locally exact from the source carpet
$K_1$ to the target carpet $K_2$.
\end{lemma}

\begin{proof}
Let $p_A$ be the fixed point of $A$. Repeating its retained word shows
that $p_A\in K_1$, and hence $F(p_A)\in K_2$ is finite. Choose a fixed
compact neighborhood of $p_A$ that avoids the pole of $F$, and then
shrink it once more. On the resulting neighborhood, $F$ has a uniform
Euclidean quasisymmetric control function $\eta$. Indeed, $F$ is a fixed
spherical quasisymmetry, and the spherical and Euclidean metrics are
bi-Lipschitz equivalent on the compact source and image neighborhoods
under consideration. This gives the required local control without
assuming that $F$ is Euclidean quasisymmetric on all of $Q$. For large
$r$, every source point used below lies in this fixed neighborhood.

Put $X=A^r(Q)$ and $d=\diam F(X)>0$, and let $x_0$ be the center of
$X$. The image $F(X)$ is a bounded closed Jordan disk, and its diameter
is attained by points on $F(\partial X)$. Thus there is
$z\in\partial X$ such that $|F(z)-F(x_0)|\ge d/2$. For every
$w\in\partial X$, the square geometry gives
$|z-x_0|/|w-x_0|\le\sqrt2$, so that
\[
 |F(w)-F(x_0)|\ge\frac{d}{2\eta(\sqrt2)}.
\]
Set
\begin{equation}\label{reset:disk-delta}
 \delta_0=\frac1{2\eta(\sqrt2)}.
\end{equation}
The Jordan curve theorem now implies
$B(F(x_0),\delta_0d)\subset F(\intr X)$.

We first choose a scale in the \emph{target} grid. Select $j\ge0$ and
write $\lambda=m^{-j}$ so that
\begin{equation}\label{reset:lambda-choice}
 \frac{\lambda}{3m}\le d<\frac\lambda3.
\end{equation}
Since $d\to0$, this is possible after discarding finitely many $r$.
Fix an integer $\ell\ge1$ such that
\begin{equation}\label{reset:ell-choice}
 \sqrt2\,m^{-\ell}<\frac{\delta_0}{3m},
 \qquad \frac13+2m^{-\ell}<1,
 \qquad \varepsilon=\lambda m^{-\ell}.
\end{equation}
Let $\mathcal G_\lambda$ be the union of all grid lines of the square
grid with mesh $\lambda$. We claim that there is
$a\in F(\partial X)$ satisfying
\begin{equation}\label{reset:far-grid}
 \dist(a,\mathcal G_\lambda)\ge\varepsilon.
\end{equation}
Otherwise the entire Jordan curve $F(\partial X)$ lies in the
$\varepsilon$-neighborhood of the grid lines. Since
$d+2\varepsilon<\lambda$, the curve meets at most one vertical grid
strip and at most one horizontal grid strip. It is therefore contained
in one strip or in the cross formed by two such strips. Every component
of the complement of that cross is unbounded. Consequently, the bounded
Jordan disk enclosed by the curve is also contained in the cross.
The distance from any point of the union of two perpendicular strips of
width $2\varepsilon$ to its complement is at most
$\sqrt2\varepsilon$. Such a cross cannot contain a disk of radius
$\delta_0d$, because \eqref{reset:lambda-choice}--\eqref{reset:ell-choice}
give $\sqrt2\varepsilon<\delta_0d$. This proves
\eqref{reset:far-grid}.

Next fix a depth in the \emph{source} grid: choose $s\ge1$ with
\begin{equation}\label{reset:s-choice}
 \eta(4n^{-s})<m^{-\ell}.
\end{equation}
Choose a relative depth-$s$ boundary cell $Y$ of $X$ containing
$F^{-1}(a)$. The full frame condition guarantees that this cell is
retained. Thus $Y=A^rC_r(Q)$ for some source word contraction $C_r$ of
length $s$. For $u,v\in Y$, choose $z\in X$ such that
$|z-u|\ge\diam X/2$. Then
$|u-v|/|u-z|\le2n^{-s}$, and a single quasisymmetric estimate gives
\[
 |F(u)-F(v)|\le\eta(2n^{-s})|F(u)-F(z)|
 \le\eta(4n^{-s})d<\frac{\varepsilon}{3}.
\]
Taking the supremum over $u,v$ yields
$\diam F(Y)<\varepsilon/3$. It follows that $F(Y)$ lies in the
interior of the unique target depth-$j$ cell $Z=B(Q)$ containing $a$,
at distance at least $2\varepsilon/3$ from its boundary. Since
$a\in K_2\cap\intr Z$, this target cell is retained.

Let $y_0,y_1$ be opposite vertices of $Y$. For every $z\in X$ we have
$|z-y_0|/|y_1-y_0|\le n^s$, whence
\[
 d\le2\eta(n^s)|F(y_1)-F(y_0)|
 \le2\eta(n^s)\diam F(Y).
\]
Together with \eqref{reset:lambda-choice}, this gives the normalized
estimates
\begin{equation}\label{reset:explicit-c0c1}
 \diam\bigl(B^{-1}F(Y)\bigr)\ge
 c_0:=\frac1{6m\eta(n^s)},\qquad
 \dist\bigl(B^{-1}F(Y),O\bigr)\ge
 c_1:=\frac23m^{-\ell}.
\end{equation}
There are only finitely many relative depth-$s$ boundary words. One of
them, denoted by $C$, occurs for infinitely many $r$. Pass to this
subsequence, write its depths as $r_i$, and denote the corresponding
target cell contractions by $B_i$. We obtain \eqref{reset:Hi} and all
the estimates in \eqref{reset:normal-bounds} concerning $H_i$.

We verify explicitly that the normal family does not degenerate. All
$H_i$ have the same spherical quasiconformal bound. Their restrictions
to $Q$ have the same Euclidean quasisymmetric control $\eta$, since the
normalizations before and after $F$ are similarities. Take the three
fixed vertices $q_0=(0,0)$, $q_1=(1,0)$, and $q_2=(0,1)$. For
$a\ne b$ and $z\in Q$, we have
$|z-q_a|/|q_b-q_a|\le\sqrt2$. Hence
\[
 |H_i(q_a)-H_i(q_b)|\ge
 \frac{\diam H_i(Q)}{2\eta(\sqrt2)}
 \ge\frac{c_0}{2\eta(\sqrt2)}.
\]
The three image points lie in the fixed compact subset
$\{z\in Q:\dist(z,O)\ge c_1\}$ and are uniformly separated.
After passing to a further subsequence, their limits exist and remain
distinct. The three-point normalized quasiconformal compactness theorem
therefore gives a nondegenerate spherical quasiconformal limit $H$ and
uniform spherical convergence of the inverse maps to $H^{-1}$. If the
maps reverse orientation, apply that theorem after composition with a
fixed reflection and then compose back. The positive boundary distance
$c_1$ implies $H(Q)\subset U$.

It remains to prove both directions of exactness. For $x\in U$, apply
the exact tail identities to the source cell $A^{r_i}C(Q)$ and to the
target cell $B_i(Q)$. Since $F(K_1)=K_2$ and $H_i(U)\subset U$, we get
\[
 x\in K_1
 \quad\Longleftrightarrow\quad
 F(A^{r_i}C(x))\in K_2
 \quad\Longleftrightarrow\quad H_i(x)\in K_2.
\]
This proves \eqref{reset:exact-Hi}. If $x\in U\cap K_1$, then
$H_i(x)\to H(x)$ and the closedness of $K_2$ imply $H(x)\in K_2$.
Conversely, let $y\in H(U)\cap K_2$ and put $x=H^{-1}(y)\in U$.
Inverse convergence gives $x_i=H_i^{-1}(y)\in U$ for all sufficiently
large $i$. By \eqref{reset:exact-Hi}, $x_i\in K_1$. Passing to the
limit and using the closedness of $K_1$, we obtain $x\in K_1$.
This proves the reverse inclusion in \eqref{reset:exact-H}.
\end{proof}

\begin{proposition}[Global contraction relations]\label{reset:power}
Under the assumptions of Lemma~\ref{reset:selection}, there are integers
$k,l\ge1$ and $b\in\Z[1/m]^2$ such that the positive homothety
\begin{equation}\label{reset:power-L}
 L(z)=m^{-l}z+b
\end{equation}
satisfies
\begin{equation}\label{reset:power-identity}
 f\circ A^k=L\circ f\quad\text{on }K_1,
 \qquad L(Q)\subset Q.
\end{equation}
\end{proposition}

\begin{proof}
Use the maps $H_i,H$ supplied by Lemma~\ref{reset:selection}, and set
$V=H(U)$. By \eqref{reset:exact-H} and
Lemma~\ref{geo:connectivity}, the set
$V\cap K_2=H(U\cap K_1)$ is connected. The full frame condition
provides the interior source contraction
\begin{equation}\label{reset:u}
 \begin{aligned}
 u&=\varphi_{(0,0)}\circ\varphi_{(n-1,n-1)},\\
 u(z)&=n^{-2}z+(n-1)n^{-2}(1,1),\\
 p_*&=(1/(n+1),1/(n+1)).
 \end{aligned}
\end{equation}
It satisfies $u(Q)\subset U$, $p_*\in U\cap K_1$, and $u(p_*)=p_*$.
On $V$ define
\[
 G_i=H_i\circ H^{-1},\qquad v=H\circ u\circ H^{-1}.
\]
Equations~\eqref{reset:exact-Hi}--\eqref{reset:exact-H} imply
\begin{align*}
 G_i(V\cap K_2)&=H_i(U)\cap K_2=G_i(V)\cap K_2,\\
 v(V\cap K_2)&=v(V)\cap K_2.
\end{align*}
Thus these are locally exact, orientation-preserving quasiconformal maps
for the same \emph{target} carpet $K_2$. Their orientation is preserved
even when $F$ or $H$ reverses orientation. Moreover, $G_i\to\id$
locally uniformly, while $v$ fixes $H(p_*)$ and is not the identity on
$V\cap K_2$. The latter assertion follows because $u$ has only one
fixed point in the plane, whereas $U\cap K_1$ contains more than one
point.

Let $\beta$ be an orientation-preserving round uniformization of the
target, followed by a M\"obius normalization for which
$S=\beta(K_2)\subset\mathbb C$. By Lemma~\ref{geo:geometry}, the
uniformization theorem, and quasisymmetric invariance of local porosity,
$S$ has measure zero and is locally porous. On $\beta(V)\cap S$, the
maps $\beta G_i\beta^{-1}$ and $\beta v\beta^{-1}$ are Schottky maps
by \cite[Theorem~4.1]{BM2020}. Each map in the first sequence is a
homeomorphism onto $\beta(H_i(U))\cap S$, where $\beta(H_i(U))$ is
open. The sequence converges locally uniformly to the identity
homeomorphism, and its relative domain is connected. The other map is
nonidentity and has the fixed point specified above. All hypotheses of
\cite[Theorem~4.2]{BM2020} are therefore satisfied. It follows that
$G_i=\id$ on $V\cap K_2$ for all sufficiently large $i$.
Consequently, $H_i=H$ on $U\cap K_1$, and then on $K_1$ by density
and continuity. Stabilization has taken place on $K_2$; the normalized
maps $H_i$ have not been regarded as self-maps of $K_1$.

Choose sufficiently large $i<j$ and put $k=r_j-r_i$ and
$L=B_jB_i^{-1}$. The equality $H_i=H_j$ gives
\[
 F A^k=L F\quad\text{on }A^{r_i}C(K_1).
\]
Hence $\Phi=F A^kF^{-1}$ and $L$ agree on a nonempty relatively open
subset of the target carpet. This subset may be chosen in
$F(\intr(A^{r_i}C(Q)))$, so it avoids the peripheral circle $f(O)$.
We have $\Phi(K_2)\subset K_2$, and $\Phi$ is locally exact on
$K_2\setminus f(O)$. Within the interior of a retained cell, the actual
holes are precisely copies of the source interior holes. It follows
that $A^k$ sends nonperipheral points to nonperipheral points, and the
same is true of $\Phi$. On the other hand, the homothety
$L=B_jB_i^{-1}$ has ratio in $m^{\Z}$ and translation vector in
$\Z[1/m]^2$. Applying Corollary~\ref{geo:onecircle} in the
\emph{target} $m$-grid gives $\Phi=L$ on all of $K_2$. This proves
the identity in \eqref{reset:power-identity}.

Since $L(K_2)=f(A^k(K_1))\subset K_2$ and $\conv K_2=Q$, we obtain
$L(Q)\subset Q$. The positive ratio of $L$ is at most $1$. If it were
$1$, then $L$ would be a translation taking $Q$ into itself, and hence
would be the identity. The relation $fA^k=f$ would then imply
$A^k(K_1)=K_1$, contradicting strict decrease of diameter. The ratio is
therefore $m^{-l}$ for an integer $l\ge1$, as required.
\end{proof}

\section{The outer square, its corners, and a plane extension}
\label{reset:outer-section}

\begin{proposition}\label{reset:outer}
Every $f\in\QS(K_1,K_2)$ satisfies $f(O)=O$ and permutes the four sides
and the four corners of $Q$. In particular,
there is a unique $R\in D_4$ agreeing with $f$ at the four corners.
Moreover, $f$ has a plane quasiconformal extension $F$ satisfying
$F(Q)=Q$ and $F(\infty)=\infty$. For such an extension,
Proposition~\ref{reset:power} applies to every source cell contraction
$A$ associated with a nonempty retained word.
\end{proposition}

\begin{proof}
We initially work with the spherical extension. The point $p_*$ in
\eqref{reset:u} belongs to $K_1$, and neither of its coordinates lies
in $\Z[1/n]$. Every side of every source peripheral polygon lies on a
grid line at some finite $n$-adic level. Thus $p_*$ is nonperipheral.
Its image $f(p_*)$ is nonperipheral in the target carpet and, in
particular, lies in $U$. Since $u^r(Q)$ shrink to $p_*$, we have
$F(u^r(Q))\subset U$ for all sufficiently large $r$. We may therefore
apply Proposition~\ref{reset:power} to obtain $fu^k=Lf$, where
$L(Q)\subset Q$ and the ratio of $L$ is strictly less than $1$.

Suppose that $f(O)\ne O$. Then $J=f^{-1}(O)$ is an interior
peripheral circle of $K_1$, and $u^k(J)$ is also an interior
peripheral circle of $K_1$. Hence
\[
 L(O)=f(u^k(J))
\]
is a peripheral circle of $K_2$. Its diameter is smaller than
$\diam O$, so it is not the outer circle and must bound an interior
hole. Since $L(Q)$ is a closed Jordan disk, that hole is exactly
$L(U)$. Choose $x\in K_2\cap U$. Then
$L(x)\in L(K_2)\subset K_2$, but also $L(x)\in L(U)$, a
contradiction. Thus $f(O)=O$. The spherical homeomorphism $F$ maps the
two components of the complement of $O$ to those components. Since
the nonempty set $K_1\cap U$ maps into $K_2\cap U$, it follows that
$F(U)=U$ and $F(Q)=Q$.

We can now choose an extension fixing infinity. If the current $F$
already fixes infinity, there is nothing to change. Otherwise let
$a=F^{-1}(\infty)\in\widehat{\mathbb C}\setminus Q$. Choose a smooth
arc from $\infty$ to $a$ in the connected open set
$\widehat{\mathbb C}\setminus Q$. A smooth time-dependent vector
field supported in a tubular neighborhood of this arc produces an
orientation-preserving smooth diffeomorphism $P$ of the sphere with
$P|_Q=\id$ and $P(\infty)=a$. A smooth diffeomorphism of the sphere
is quasiconformal. Replacing $F$ by $F\circ P$ leaves it unchanged on
$Q$ and makes it fix infinity. We henceforth use this plane
quasiconformal extension, still denoted by $F$. It satisfies $F(Q)=Q$
and $F|_{K_1}=f$. Every source cell $A^r(Q)$ is contained in $Q$,
so the hypothesis of Proposition~\ref{reset:power} now holds for every
source cell contraction.

We next identify the sides. Consider the source bottom side
$e=[0,1]\times\{0\}$. For $d\ge1$ and $1\le a\le n^d-2$, the full
frame condition guarantees that
\[
 A_{d,a}(z)=n^{-d}z+(a/n^d,0)
\]
is a retained source cell contraction. Its fixed point is
$p_{d,a}=(a/(n^d-1),0)$, and these points are dense in $e$. Choose
one such point $p$ whose image $q=f(p)$ is not a target corner.
Since $f$ is injective, at most four choices of $p$ are excluded.
Let $e'$ be the target side whose relative interior contains $q$.
The contraction relation gives $fA_{d,a}^k=Lf$, where $L$ has ratio
$s\in(0,1)$ and fixes $q$. Thus $L(z)=sz+(1-s)q$.
If $q$ is in the interior of the bottom side, the horizontal projection
of $L(Q)$ is contained strictly in $(0,1)$ and its vertical projection
is $[0,s]$. Consequently, $L^{-1}(O)\cap Q$ is exactly the bottom
side. The other three cases are identical after a change of coordinates.
In general,
\[
 L^{-1}(O)\cap Q=e'.
\]
For $x\in e$, we have $A_{d,a}^k(x)\in O$, and therefore
$L(f(x))\in O$. It follows that $f(e)\subset e'$.
Applying the same argument to all four source sides, we see that their
images cover $O$ and that each image lies in a single target side.
The interior of every target side must be covered. Thus the four sides
correspond bijectively, and each maps onto the corresponding target
side. Their common endpoints map to the corresponding corners.
Finally, a homeomorphism of $O$ preserves or reverses cyclic order, so
this corner permutation determines a unique element $R\in D_4$.
\end{proof}

\section{Translations in a thin strip and linearity on the sides}
\label{reset:boundary-section}

\begin{proposition}\label{reset:boundary}
Let $R(z)=A_Rz+c_R\in D_4$ be the corner symmetry supplied by
Proposition~\ref{reset:outer}. Then every $f\in\QS(K_1,K_2)$ satisfies
\begin{equation}\label{reset:boundary-identity}
 f|_O=R|_O.
\end{equation}
\end{proposition}

\begin{proof}
Use the plane extension $F$ chosen in Proposition~\ref{reset:outer}.
For each integer $s\ge1$, put $a=n^{-s}$ and consider the two source
bottom-row contractions
\[
 D_a(z)=az,\qquad C_a(z)=az+ae_1,\qquad e_1=(1,0).
\]
Proposition~\ref{reset:power} gives a contraction relation for a
positive power of each map. Take a common multiple $B\ge1$ of those
two exponents and iterate the corresponding target homotheties. This
gives
\begin{equation}\label{reset:two-powers}
 fD_a^B=L_Df,\qquad fC_a^B=L_Cf,
\end{equation}
where $L_D,L_C$ have ratios in $m^{-\N}$ and translation vectors in
$\Z[1/m]^2$. No bound on $B$ independent of $s$ is required. Set
\begin{equation}\label{reset:thin-parameters}
 \delta=a^B,\qquad t=a+a^2+\cdots+a^B,\qquad T(z)=z+te_1.
\end{equation}
Then $C_a^B=TD_a^B$, and
\[
 \delta^{-1},\ t/\delta,\ (1-t)/\delta\in\N,
 \qquad t+\delta\le\frac a{1-a}+a\le\frac56<1.
\]
By \eqref{reset:two-powers}, the orientation-preserving plane
quasiconformal map $\Phi=FTF^{-1}$ agrees with the positive homothety
$L=L_CL_D^{-1}$ on $F(D_a^B(K_1))$. Here
\begin{equation}\label{reset:thin-target-arithmetic}
 \operatorname{ratio}(L)\in m^{\Z},\qquad
 \operatorname{trans}(L)\in\Z[1/m]^2.
\end{equation}

Consider the open strip
\[
 W_s=(0,1-t)\times(-\delta,\delta),\qquad V_s=F(W_s).
\]
Within $0\le y<\delta$, every bottom-row cell of side length
$\delta$ is retained, and the pattern in each cell is a copy of
$K_1$ with the same orientation. Since $t/\delta$ is an integer,
$T$ translates each such pattern exactly onto the pattern of another
bottom-row cell. The patterns also match on both sides of every common
vertical edge. Below $y=0$, both relevant neighborhoods are disjoint
from $K_1$. The upper line $y=\delta$ and the two vertical end lines
are excluded from $W_s$, so no assertion about the next row or an outer
corner is involved. Thus $T$ is locally exact at every point of
$W_s\cap K_1$, and $T(W_s\cap K_1)\subset K_1$.

We also verify connectedness of this relative domain. Let
$q=(1-t)/\delta\in\N$, and form the finite row
\[
 Z_s=\bigcup_{j=0}^{q-1}\bigl(\delta K_1+j\delta e_1\bigr).
\]
In each cell, the carpet with its own outer boundary removed is
connected and dense, by Lemma~\ref{geo:connectivity}. Adjoining the
common open vertical edges joins these sets into the connected set
obtained from $Z_s$ by deleting its rectangular outer boundary.
Adjoining the interior points of the bottom side preserves
connectedness, since all these points lie in the closure of that set.
The resulting set is precisely $W_s\cap K_1$. Hence
$V_s\cap K_2=F(W_s\cap K_1)$ is connected.

It follows that $\Phi$ is locally exact at every point of
$V_s\cap K_2$ and maps this set into $K_2$. Since
$D_a^B(U)\subset W_s$, the initial equality with $L$ holds on a
nonempty relatively open subset of this connected relative domain.
Apply Lemma~\ref{geo:propagation} using the \emph{target} $m$-grid;
the alignment required there follows from
\eqref{reset:thin-target-arithmetic}. We obtain $\Phi=L$ throughout
$V_s\cap K_2$. In particular, continuity at the two ends of the bottom
interval gives
\begin{equation}\label{reset:boundary-increment}
 f(x+t,0)=L(f(x,0)),\qquad 0\le x\le1-t.
\end{equation}

The sides and corners have already been identified. We may therefore
write, uniquely,
\[
 f(x,0)=c_R+(A_Re_1)h(x),\qquad 0\le x\le1,
\]
where $h:[0,1]\to[0,1]$ is an increasing homeomorphism with
$h(0)=0$ and $h(1)=1$. Taking $x=0$ in
\eqref{reset:boundary-increment} eliminates the translation vector of
$L$. Denote its ratio by $\lambda_s$ and the translation length in
\eqref{reset:thin-parameters} by $t_s$. Then
\begin{equation}\label{reset:interval-affine}
 h(x+t_s)=\lambda_s h(x)+h(t_s),\quad
 0\le x\le1-t_s,\qquad \lambda_s\in m^{\Z}.
\end{equation}
The source parameter $n$ produces the small translations, whereas the
discrete ratios come from the target parameter $m$. Equality of these
parameters has not been assumed. Regardless of the dependence of $B$
on $s$, we have
\[
 n^{-s}\le t_s\le\frac{n^{-s}}{1-n^{-s}}\longrightarrow0.
\]
For large $s$, substitute $x=1/2$ into
\eqref{reset:interval-affine}. Continuity gives
\[
 \lambda_s=\frac{h(1/2+t_s)-h(t_s)}{h(1/2)}\longrightarrow1.
\]
Since $m^{\Z}$ is discrete near $1$, we have $\lambda_s=1$ for all
sufficiently large $s$. Thus
\begin{equation}\label{reset:interval-additive}
 h(x+t_s)=h(x)+h(t_s),\qquad 0\le x\le1-t_s.
\end{equation}

For completeness, the following bounds determine $h$ without any
differentiability assumption. Put $q_s=\lfloor1/t_s\rfloor$.
Iterating \eqref{reset:interval-additive}, we obtain
$h(jt_s)=jh(t_s)$ for $0\le j\le q_s$, and hence
$q_sh(t_s)\le1$. Also, $1-t_s\le q_st_s$, so evaluation at
$x=1-t_s$ yields
\[
 1=h(1-t_s)+h(t_s)\le(q_s+1)h(t_s).
\]
Since $q_st_s\to1$, it follows that
\[
 \frac1{(q_s+1)t_s}\le\frac{h(t_s)}{t_s}
 \le\frac1{q_st_s},\qquad \frac{h(t_s)}{t_s}\longrightarrow1.
\]
For a fixed $0\le x<1$, put $p_s=\lfloor x/t_s\rfloor$. For all
large $s$, $(p_s+1)t_s\le1$, and therefore
\[
 p_sh(t_s)\le h(x)\le(p_s+1)h(t_s).
\]
Both bounds tend to $x$, so $h(x)=x$. The endpoint $x=1$ is already
fixed by normalization. Thus $f$ agrees with $R$ on the bottom side.

Apply square isometries separately to the source and target coordinates
to send any corresponding pair of sides to the bottom side, and repeat
the argument. The transformed digit sets still define fixed patterns
with the same orientation at every stage, with subdivision parameters
$n$ and $m$, respectively, and they still satisfy
condition~\eqref{main:P}. No symmetry of either original pattern is
needed for this change of coordinates. We conclude that the parameter
on each side is the linear one prescribed by its corner images, and
hence $f|_O=R|_O$.
\end{proof}

\begin{figure}[htbp]
\centering
\begin{tikzpicture}[scale=7.4]
 \fill[blue!7] (0,0) rectangle (.7,.1);
 \draw[gray!55,thin,step=.1] (0,0) grid (1,.1);
 \draw[thick] (0,0) rectangle (1,.48);
 \draw[dashed] (0,.1)--(1,.1);
 \fill[blue!22] (0,0) rectangle (.1,.1);
 \fill[orange!25] (.3,0) rectangle (.4,.1);
 \draw[thick,blue!70!black] (0,0) rectangle (.1,.1);
 \draw[thick,orange!80!black] (.3,0) rectangle (.4,.1);
 \draw[-{Stealth[length=2mm]},thick] (.05,.28)--node[above]{$T:z\mapsto z+te_1$}(.35,.28);
 \draw[<->] (0,-.065)--node[below]{$1-t$}(.7,-.065);
 \draw[<->] (1.035,0)--node[right]{$\delta$}(1.035,.1);
 \node[above,font=\small] at (.08,.12){$D_a^B Q$};
 \node[above,font=\small] at (.35,.12){$C_a^B Q$};
 \node at (.65,.40){Boundary strip};
 \node[below left] at (0,0){$0$};
 \node[below right] at (1,0){$1$};
\end{tikzpicture}
\caption{Translation along the bottom row of retained cells.
The proportions are schematic; the parameters satisfy
\eqref{reset:thin-parameters}. The top edge of the strip is excluded
from the local exactness domain. The boundary identity is obtained by
continuity.}
\label{unified:stripfigure}
\end{figure}
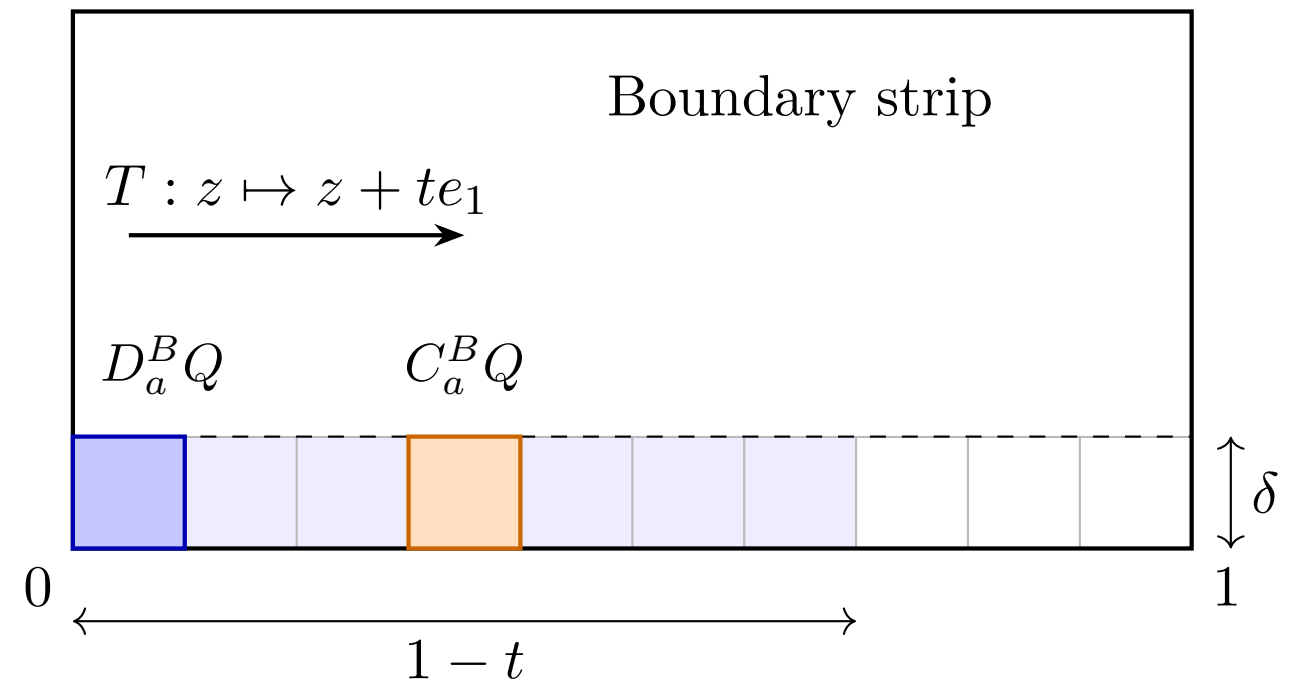

\section{Translation-periodic extension and global rigidity}

We do not yet assume that $n$ and $m$ have a common power.
Let $f:K_1\to K_2$ be a quasisymmetric homeomorphism.
By Proposition~\ref{reset:boundary}, there is a unique
\[
 R(z)=Az+c\in D_4
\]
such that $f|_O=R|_O$, where $A$ is a signed permutation matrix and
$c\in\{0,1\}^2$.
We have also chosen a quasiconformal extension $F$ with
$F|_{K_1}=f$ and $F(Q)=Q$.
Its pole lies outside a neighborhood of $Q$.

\subsection{Extension by integer translations}

Define
\begin{equation}\label{periodic:sets}
 \mathcal K_i=\bigcup_{p\in\Z^2}(K_i+p),\qquad i=1,2.
\end{equation}
All copies in this definition have the same orientation.
We use translations rather than reflections, since a general pattern
need not be invariant under reflection.

For $x\in Q$ and $p\in\Z^2$, set
\begin{equation}\label{periodic:extension}
 \widehat F(x+p)=F(x)+Ap.
\end{equation}
This definition agrees on common sides. Indeed, if $x,y\in O$ and
$x+p=y+q$, then
\[
 F(x)+Ap=Ax+c+Ap=Ay+c+Aq=F(y)+Aq.
\]
Each closed unit square $Q+p$ is mapped homeomorphically onto $Q+Ap$.
The target squares tile the plane, and their side and vertex
identifications agree. Thus $\widehat F$ is a plane homeomorphism.
Its quasiconformal dilatation is bounded by the same constant in every
open unit square.
To verify quasiconformality across the integer grid, work in a compact
rectangle. Almost every horizontal or vertical segment meets only
finitely many seams. Absolute continuity on its constituent segments,
together with continuity at the seams, gives absolute continuity on the
whole segment. Hence $\widehat F$ is ACL and satisfies a uniform
quasiconformal distortion bound off the integer grid, which has area
zero. It is therefore quasiconformal on the plane; the
orientation-reversing case is treated in the same way.
By construction, it maps $\mathcal K_1$ onto $\mathcal K_2$.

Set
\begin{equation}\label{periodic:G}
 G(z)=\widehat F(z)-c.
\end{equation}
Since $c\in\Z^2$, this translation preserves the target periodic set.
Consequently,
\begin{equation}\label{periodic:translation}
 \begin{gathered}
 G(\mathcal K_1)=\mathcal K_2,\qquad G(0)=0,\\
 G(z+p)=G(z)+Ap\quad(z\in\C,\ p\in\Z^2),\\
 M:=\sup_{z\in\C}|G(z)-Az|<\infty.
 \end{gathered}
\end{equation}
The last assertion follows because $G(z)-Az$ is integer-periodic and
continuous on $Q$.

\subsection{Local exactness of expanding maps}

For the moment, let $K=K(q,\mathcal F)$ be either fixed-pattern carpet,
and write
$\mathcal K=\bigcup_{p\in\Z^2}(K+p)$.
For $r\ge1$, let $T_{q,r}(z)=q^rz$.

\begin{lemma}[Expansion and finite-generation exceptions]\label{periodic:expansion}
We have $T_{q,r}(\mathcal K)\subset\mathcal K$.
Let $\mathcal E_{q,r}$ consist of the boundaries of all holes of
generations $1,\ldots,r$ in every integer unit square.
This family is locally finite, and $T_{q,r}$ is locally exact for
$\mathcal K$ at every point of
$\mathcal K\setminus\bigcup_{J\in\mathcal E_{q,r}}J$.
\end{lemma}

\begin{proof}
Given $x\in K$, choose a retained generation-$r$ cell containing $x$.
The full frame condition and the exact-copy identity in
Lemma~\ref{geo:geometry} give $q^rx\in K+d$, where
$d\in\Z^2$ is the digit of this cell.
For $x+p\in K+p$, adding the integer vector $q^rp$ proves the forward
inclusion.

At a point in the interior of a retained generation-$r$ cell, local
exactness is precisely the exact-copy identity for its tail.
Suppose instead that the point lies on a generation-$r$ grid side or
vertex, but not on any hole boundary of the first $r$ generations.
The two or four adjacent generation-$r$ cells are then all retained.
Their tails are copies of $K$ with the same orientation. After expansion,
they agree exactly with the periodic pattern on the two sides of an
integer grid edge or in the four quadrants at an integer grid vertex.
The same argument applies on the sides and vertices of the original unit
squares, since the periodic construction supplies a complete unit
square of the same orientation on each side.
This proves both directions of local exactness at such points.
There are only finitely many holes of the first $r$ generations in each
unit square, and their closures lie strictly inside that square.
Thus $\mathcal E_{q,r}$ is locally finite.
\end{proof}

Only forward invariance under expansion is used here.
The inclusion $q^{-r}\mathcal K\subset\mathcal K$ need not hold and
is not part of the argument.

\subsection{Uniqueness on the periodic carpet}

To apply Schottky-map uniqueness on an unbounded set, we first justify
round uniformization and the connectedness of the domain obtained by
removing the exceptional hole boundaries. In particular, the finite
peripheral-circle removal lemma must not be applied directly to an
arbitrary infinite family.

\begin{lemma}[The periodic analytic domain]\label{periodic:analysisdomain}
For a fixed pattern satisfying~\eqref{main:P}, the compactification
$X=\mathcal K\cup\{\infty\}$ is a spherical Sierpi\'nski carpet of
area zero. Its peripheral circles are uniform quasicircles and are
uniformly relatively separated in the spherical metric.
There is a plane quasiconformal homeomorphism $\beta$ such that
$S=\beta(\mathcal K)$ is a relative Schottky set in $\C$ of area zero
that is locally porous at every finite point.
No porosity assumption at $\infty$ is required.

Let $\mathcal E$ be a locally finite family of boundaries of actual
holes of $\mathcal K$, and let $V_J$ be the hole bounded by $J$.
Set
\begin{equation}\label{periodic:W}
 W=\C\setminus\bigcup_{J\in\mathcal E}\overline{V_J}.
\end{equation}
Then $W$ is open, and
\begin{equation}\label{periodic:connected}
 W\cap\mathcal K
   =\mathcal K\setminus\bigcup_{J\in\mathcal E}J
\end{equation}
is connected and dense in $\mathcal K$.
\end{lemma}

\begin{proof}
The complementary components of the periodic set are exactly the
interior holes in its integer unit squares. They are positive homothetic
copies and integer translates of finitely many Jordan polygons, and
their closures are pairwise disjoint.
Every hole has Euclidean diameter at most $D_0=\sqrt2$.
Relative separation for holes in one unit square follows from
Lemma~\ref{geo:geometry}.
For holes in different unit squares, take the one with the smaller
parent cell. Its distance from the boundary of that cell is at least a
fixed positive multiple of the cell side length, whereas the other hole
lies outside the interior of that cell.
The entire periodic hole family is therefore uniformly relatively
separated in the Euclidean metric.

A uniform metric comparison gives relative separation on the sphere.
For two closed holes $B_1,B_2$, put
$a_i=(1+\dist(0,B_i)^2)^{1/2}$ and assume $a_1\le a_2$.
The spherical length factor on $B_i$ is comparable to $a_i^{-2}$,
with constants depending only on $D_0$.
If $a_2/a_1$ is bounded by a fixed constant, Euclidean relative
separation gives spherical relative separation.
If $a_2>2(1+D_0)a_1$, the Euclidean distance between the holes is at
least $a_2/4$. Their spherical distance is therefore at least
$C^{-1}a_1^{-1}$, whereas the spherical diameter of the farther hole is
at most $D_0a_2^{-2}$. The ratio again has a uniform positive lower bound.
These estimates follow from the chordal-distance formula
\[
 \chi(z,w)=\frac{|z-w|}{\sqrt{1+|z|^2}\sqrt{1+|w|^2}},
\]
and $C$ depends only on $D_0$.
The uniform quasicircle constant comes from the finite collection of
root polygons. Translations and positive homotheties are M\"obius maps
of the sphere and preserve the corresponding quasicircle constants.

The spherical diameters of the holes form a null family.
In a bounded region, their diameters tend to zero with their generation;
for holes escaping to infinity, the uniform Euclidean diameter bound
gives the same conclusion.
A finite integer rectangle filled with copies of $K$ of the same
orientation is an ordinary rectangular carpet. Thus $X$ is locally
connected at every finite point.
At infinity, use relative open neighborhoods outside integer squares.
The edge-adjacency graph of the unit squares outside a large square is
connected, and the full frames join the carpet copies in those squares.
These neighborhoods are therefore connected, proving local
connectedness at infinity.
Together with the Jordan complementary components, the null-family
property, and the density of the holes, this shows that $X$ is a
spherical Sierpi\'nski carpet.
It has area zero because $\mathcal K$ is a countable union of sets of
area zero.

Bonk's round uniformization theorem
\cite[Theorem~1.1, Corollary~1.2]{Bonk2011} maps $X$ to a round carpet
by a spherical quasiconformal homeomorphism.
Postcompose with a M\"obius map taking the image of infinity to infinity.
The resulting map $\beta$ is a plane quasiconformal homeomorphism.
Since the original point at infinity belongs to no peripheral circle,
the round complementary components are bounded disks.
Thus $S=\beta(\mathcal K)$ is a relative Schottky set in $\C$.
It has area zero. The porosity estimate of Lemma~\ref{geo:geometry},
applied in the integer squares, gives local porosity of $\mathcal K$
at every finite point. A plane quasiconformal map is locally
quasisymmetric, so $S$ is locally porous as well.
No porosity at the point at infinity has been used.

It remains to prove~\eqref{periodic:connected}.
Local finiteness of $\mathcal E$ implies local finiteness of the
associated closed holes: the holes have uniformly bounded diameters,
so if one meets a compact set, its boundary meets the $D_0$-neighborhood
of that set. Hence~\eqref{periodic:W} defines an open set.
For $M\ge1$, let $B_M=[-M,M]^2$ and put
\[
 C_M=(\mathcal K\cap\intr B_M)
       \setminus\bigcup_{J\in\mathcal E}J.
\]
Every hole lies strictly inside an integer unit square, so none crosses
$\partial B_M$. Only finitely many members of $\mathcal E$ occur in
$B_M$.
Apply Lemma~\ref{geo:connectivity} to the rectangular carpet
$\mathcal K\cap B_M$, removing those finitely many peripheral circles
and the outer boundary of the rectangle.
It follows that $C_M$ is connected and dense in that rectangular carpet.
The sets $C_M$ increase with $M$, and
\[
 \bigcup_{M\ge1}C_M=W\cap\mathcal K.
\]
Their union is connected and dense in $\mathcal K$.
Thus the finite-circle removal lemma is applied only in finite
rectangles; the unbounded conclusion follows by an increasing union.
\end{proof}

\begin{lemma}[Continuation on the periodic carpet]\label{periodic:uniqueness}
Let $H_1,H_2$ be orientation-preserving quasiconformal homeomorphisms of
the plane, both mapping $\mathcal K$ into itself.
Suppose that both are locally exact outside a locally finite family
$\mathcal E$ of actual hole boundaries.
If they agree on a nonempty relatively open subset of $\mathcal K$,
then $H_1=H_2$ on $\mathcal K$.
\end{lemma}

\begin{proof}
Take $W$ and $\beta$ from Lemma~\ref{periodic:analysisdomain}.
The set $W\cap\mathcal K$ is connected and dense, so it meets the
initial relatively open equality set in a nonempty relatively open set.
Both $\beta H_i\beta^{-1}$ restrict to Schottky maps on
$\beta(W)\cap S$: they preserve orientation and are locally exact
quasiconformal homeomorphisms near every point of this relative domain.
Proposition~\ref{geo:analysis}, namely the local regularity and
uniqueness statements of
\cite[Theorems~4.1, 4.3]{BM2020}, gives equality throughout that
connected relative domain.
Density and continuity extend the equality to all of $\mathcal K$.
\end{proof}

\subsection{Comparing expansion with integer translation}

\begin{proposition}[Global isometry and a common power]\label{periodic:global}
Under the hypotheses of Theorem~\ref{main:classification}, every
$f\in\QS(K_1,K_2)$ is the restriction of an $R\in D_4$ to $K_1$.
Moreover, $n^k=m^\ell$ for some $k,\ell\ge1$.
\end{proposition}

\begin{proof}
Apply Proposition~\ref{reset:power} to the source corner contraction
$a_n(z)=z/n$.
There are $k,\ell\ge1$ and a positive homothety $L$ such that
\[
 f a_n^k=L f,\qquad \operatorname{ratio}(L)=m^{-\ell}.
\]
Since $a_n(0)=0$ and $f(0)=c$, the homothety fixes $c$; hence
$L(z)=c+m^{-\ell}(z-c)$.
On $K_1$ we have $G=f-c$, and therefore
\begin{equation}\label{periodic:localpower}
 G(n^{-k}x)=m^{-\ell}G(x),\qquad x\in K_1.
\end{equation}

Consider the plane maps
\begin{equation}\label{periodic:twoH}
 H_1=G^{-1}T_{m,\ell}G,\qquad H_2=T_{n,k}.
\end{equation}
Both preserve orientation, even when $G$ reverses orientation.
Lemma~\ref{periodic:expansion} and
$G(\mathcal K_1)=\mathcal K_2$ show that both map $\mathcal K_1$ into
itself.
The exceptional circles for $H_2$ are the source hole boundaries of the
first $k$ generations in each unit square.
For $H_1$, take the inverse images under $G$ of the target hole
boundaries of the first $\ell$ generations in each unit square.
They are actual source hole boundaries, since the plane homeomorphism
$G$ bijects the complementary components of the two periodic carpets.
This family is locally finite because $G$ is proper.
The union of the two exceptional families is therefore locally finite,
and both maps are locally exact off that union.

Equation~\eqref{periodic:localpower} gives $H_1=H_2$ on $n^{-k}K_1$.
This set contains the nonempty relatively open set
$\mathcal K_1\cap\intr(n^{-k}Q)$.
Lemma~\ref{periodic:uniqueness} extends the equality to
$\mathcal K_1$, giving
\begin{equation}\label{periodic:globalpower}
 G(n^kz)=m^\ell G(z),\qquad z\in\mathcal K_1.
\end{equation}
The comparison takes place on the source periodic carpet.
We have not assumed that multiplication by $n^k$ preserves the target
periodic carpet.

Apply~\eqref{periodic:globalpower} at $z$ and $z+p$, with
$p\in\Z^2$, and subtract.
Using~\eqref{periodic:translation}, we obtain
\[
 n^kAp=m^\ell Ap.
\]
Taking $p\ne0$ and using the invertibility of $A$ yields
\begin{equation}\label{periodic:commonpower}
 N_0:=n^k=m^\ell>1.
\end{equation}
Thus the common power is a consequence of periodicity, rather than a
hypothesis in the preceding arguments on the two grids.

Iterating~\eqref{periodic:globalpower}, for $z\in\mathcal K_1$ and
$r\ge1$ we have
\[
 G(z)-Az=N_0^{-r}
       \bigl(G(N_0^rz)-AN_0^rz\bigr).
\]
The expression in parentheses has absolute value at most the constant
$M$ in~\eqref{periodic:translation}.
Letting $r\to\infty$ gives $G(z)=Az$.
Consequently, $f(z)=Az+c=R(z)$ on $K_1$.
Since $f$ is onto, $R(K_1)=K_2$.
Conversely, the restriction of any $R\in D_4$ with $R(K_1)=K_2$ is
quasisymmetric.
This proves the mapping classification and the common-power necessity
in Theorem~\ref{main:classification}.
\end{proof}

\begin{proposition}[Opposite-side compatibility]\label{unified:opposite}
Suppose $K_1,K_2$ satisfy the geometric hypotheses of
Theorem~\ref{main:classification}.
Let $f\in\QS(K_1,K_2)$ preserve the outer boundary and agree at the four
corners with $R(z)=Az+c\in D_4$.
If
\begin{equation}\label{unified:opposite-identities}
 \begin{aligned}
 f(1,t)-f(0,t)&=Ae_1,\\
 f(t,1)-f(t,0)&=Ae_2\qquad(0\le t\le1),
 \end{aligned}
\end{equation}
then these boundary compatibility identities, together with the power
conjugacy of Proposition~\ref{reset:power}, suffice for the
periodization argument to give $f=R|_{K_1}$ and $n^k=m^\ell$ for some
$k,\ell\ge1$.
In particular, taking $K_1=K_2$ gives the opposite-side compatibility
criterion for self-maps.
\end{proposition}

\begin{proof}
Take a spherical quasiconformal extension $F$ of $f$.
Since $f(O)=O$ and carpet points inside the outer square map to carpet
points inside it, the extension maps $Q$ onto $Q$.
Identities~\eqref{unified:opposite-identities} ensure that
$\widehat F(x+p)=F(x)+Ap$ agrees on the common sides of adjacent unit
squares. The corner condition ensures agreement at vertices.
Thus~\eqref{periodic:extension}--\eqref{periodic:translation} remain
valid, without requiring $f=R$ on the whole outer boundary beforehand.
Power conjugacy for the source corner contraction and $f(0)=c$ then
give~\eqref{periodic:localpower}.
The rest of the proof of Proposition~\ref{periodic:global} uses only
these identities: continuation gives the global expansion relation,
comparison with integer translations yields $n^k=m^\ell$, and the bounded
periodic error then vanishes. Hence $f=R|_{K_1}$.
\end{proof}

\section{A finite classification criterion and purely periodic constructions}

\subsection{Recovering the digits from the limit set}

\begin{lemma}[Uniqueness of the digit set on a fixed grid]\label{periodic:digits}
Let $K=K(N,\mathcal F)$ satisfy the full frame condition.
For every $d\in\{0,\ldots,N-1\}^2$,
\begin{equation}\label{periodic:recover}
 d\in\mathcal F
 \quad\Longleftrightarrow\quad
 K\cap\intr\!\left(\frac{d+Q}{N}\right)\ne\varnothing.
\end{equation}
In particular, $K$ determines its digit set on the prescribed grid.
\end{lemma}

\begin{proof}
If the digit is not retained, the interior of its cell is disjoint from
every retained first-generation closed cell, and hence from $K$.
Conversely, the full frame condition makes
\[
 u_N(z)=N^{-2}z+(N-1)N^{-2}(1,1)
\]
a retained-cell contraction. Its fixed point
$p_*=(1/(N+1),1/(N+1))$ belongs to $K\cap U$.
If $d$ is retained, then $(d+p_*)/N$ belongs to $K$ in the interior of
that cell, proving~\eqref{periodic:recover}.
\end{proof}

\begin{proposition}[A terminating finite criterion]\label{periodic:finite}
Let $K_1,K_2$ satisfy the hypotheses of
Theorem~\ref{main:classification}.
If $n$ and $m$ have no common positive integer power, the carpets are
not quasisymmetrically equivalent.
If $N=n^a=m^b$, form the grouped digit sets in
\eqref{main:digitsD} and~\eqref{main:digitsE}.
For each $R(z)=Az+c\in D_4$,
\begin{equation}\label{periodic:finitecriterion}
 R(K_1)=K_2
 \quad\Longleftrightarrow\quad
 E^{[b]}=\{Ad+(N-1)c:d\in D^{[a]}\}.
\end{equation}
The existence of a common power and all eight digit comparisons can
be decided by finite computations.
\end{proposition}

\begin{proof}
Composing $a$ and $b$ consecutive generations, respectively, gives
\[
 K(n,D)=K(n^a,D^{[a]}),\qquad
 K(m,E)=K(m^b,E^{[b]}).
\]
The grouped patterns still satisfy the full frame condition.
Their first-generation holes are exactly the holes of the corresponding
original patterns in the first $a$ or $b$ generations.
Their closures remain pairwise disjoint, including holes of different
generations: a later hole lies strictly inside a retained cell, while
an earlier hole is disjoint from the interior of that cell.
Thus the grouped patterns also satisfy~\eqref{main:P}.
For every $d$ and $x\in Q$, direct calculation gives
\[
 R\left(\frac{d+x}{N}\right)
   =\frac{Ad+(N-1)c+R(x)}{N}.
\]
Hence $R$ permutes the $N$-grid cells by
$d\mapsto Ad+(N-1)c$.
Lemma~\ref{periodic:digits} says that equality of the resulting limit
sets is equivalent to equality of their digit sets, proving
\eqref{periodic:finitecriterion}.
Together with Proposition~\ref{periodic:global}, this proves the finite
classification assertion of Theorem~\ref{main:classification}.

To make termination explicit, write
\[
 n=\prod_p p^{\alpha_p},\qquad
 m=\prod_p p^{\beta_p}.
\]
A common power exists precisely when the supports of the two exponent
vectors agree and the vectors are positive rational multiples of one
another.
First compare their supports. If they agree, choose any prime $p_0$ in
that support and set
\[
 g=\gcd(\alpha_{p_0},\beta_{p_0}),\qquad
 a=\beta_{p_0}/g,\qquad b=\alpha_{p_0}/g.
\]
Check the finitely many equalities $a\alpha_p=b\beta_p$.
If all hold, the coprime pair $(a,b)$ gives the least positive
common-power exponents; otherwise, no common power exists.
For the resulting finite integer $N=n^a=m^b$, construct the grouped
digit sets and compare their eight $D_4$ images.
No unbounded search over exponents is needed.
\end{proof}

Once two carpets have been identified by a square symmetry, their
common-power condition also follows from the multidimensional
Cobham theorem of Chan and Hare~\cite{CH14} and Charlier, Leroy and
Rigo~\cite{CLR15}. For a fixed grid of base $n$, the sets
$(n^k K-d)\cap Q$, where $k\ge0$ and $d\in\{0,\ldots,n^k-1\}^2$,
take only finitely many values. A retained cell gives $K$ itself.
In any other cell only whole sides or vertices can occur, since all
retained cell boundaries belong to $K$. If the same carpet had two
multiplicatively independent grid bases, the cited theorem would make
it a finite union of rational polyhedra. Such a union has only finitely
many complementary components, whereas $K$ has holes in every generation.

\begin{remark}[Different grids: the parameters $3$ and $9$]
Different subdivision parameters need not describe different isometry
classes. Let
$D_3=\{0,1,2\}^2\setminus\{(1,1)\}$.
Then
\[
 K(3,D_3)=K(9,D_3^{[2]}),\qquad |D_3^{[2]}|=8^2=64.
\]
This is the standard middle-third carpet with two generations grouped
into one.
In contrast, deleting only the central cell of the $9$-grid gives
$E_9=\{0,\ldots,8\}^2\setminus\{(4,4)\}$, which has $80$ digits.
On the common $9$-grid, no square symmetry can identify a $64$-digit
pattern with an $80$-digit pattern.
Thus $K(3,D_3)$ and $K(9,E_9)$ are not quasisymmetrically equivalent.
The existence of a common power alone does not decide equivalence; the
grouped patterns must still be compared.
\end{remark}

\subsection{Constructions periodic from the first generation}

\begin{corollary}[Purely periodic classification]\label{periodic:pureperiodic}
Consider two grid constructions whose subdivision parameters and
deletion patterns are periodic from the first generation.
At each generation the subdivision parameter is at least $3$, all
retained parent cells use the same pattern with the same orientation,
and the full frame is retained.
Suppose each period contains a nonempty deletion and both limit sets
are connected and have no local cut points.
Then every quasisymmetric homeomorphism between the two limit sets is
the restriction of an element of $D_4$.
After grouping one full period in each construction, equivalence is
decided by the common-power and digit-set tests of
Proposition~\ref{periodic:finite}.
\end{corollary}

\begin{proof}
For one period, let the subdivision parameters be
$q_1,\ldots,q_t$ and the digit sets be
$\mathcal F_1,\ldots,\mathcal F_t$.
Set
\[
 N=\prod_{j=1}^tq_j,\qquad
 \mathcal F=\left\{
  \sum_{j=1}^t\left(\prod_{r=j+1}^tq_r\right)d_j:
       d_j\in\mathcal F_j\right\},
\]
where an empty product equals $1$.
Composing the contractions over one period gives the fixed-pattern
IFS $K=K(N,\mathcal F)$, with the same orientation in every retained
cell.
The full frame condition at every generation gives the full frame
condition for the grouped pattern.
A nonempty deletion during the period makes $\mathcal F$ a proper
subset of the full $N$-grid digit set.
Since the limit set is connected and has no local cut points,
Proposition~\ref{geo:nolocalcut} implies that the grouped pattern
satisfies~\eqref{main:P}.
Group the two constructions in this way and apply
Theorem~\ref{main:classification} and
Proposition~\ref{periodic:finite}.
\end{proof}

\section{Self-map groups and examples}\label{unified:applications}

The classification identifies a self-map group by a finite computation
on the first approximation. One may use either the retained digits or
the actual complementary domains.

\begin{proposition}[First-generation stabilizers]\label{unified:stabilizer}
Let $K=K(n,D)$ satisfy the hypotheses of
Theorem~\ref{main:classification}, and let
$\mathcal V=\{V_1,\ldots,V_h\}$ be its first-generation holes. Then
\begin{align}
 \QS(K)=\BiLip(K)=\Iso(K)
 &=\{R|_K:R\in D_4,\ \rho_n^R(D)=D\}\notag\\
 &=\{R|_K:R\in D_4,\ R(\mathcal V)=\mathcal V\}.
 \label{unified:stabilizerformula}
\end{align}
In particular, eight symmetry comparisons determine the group.
\end{proposition}

\begin{proof}
The first equality with a stabilizer follows from
Corollary~\ref{main:isometricclassification} and the fixed-grid digit
criterion. A square symmetry permuting the first-generation holes
preserves $Q_1$, hence permutes the retained digits and preserves $K$.
For the converse, suppose that $R(K)=K$. It permutes the complementary
domains and preserves $O$. Every first-generation hole contains an
open grid cell, so its diameter is at least $\sqrt2/n$. The closure of
every later hole lies strictly inside a parent cell of side length at
most $1/n$, and therefore has diameter strictly less than $\sqrt2/n$.
Since an isometry preserves diameters, $R$ permutes the first-generation
holes. This gives the second stabilizer description.
\end{proof}

\subsection{Separated unit-cell holes}

Suppose that the deleted digits form a nonempty set
\[
 H\subset\{1,\ldots,n-2\}^2,\qquad
 \|h-h'\|_\infty\ge2\quad(h\ne h').
\]
Each actual first-generation hole is then one open cell, and
condition~\eqref{main:P} holds. All quasisymmetric maps and self-map
groups for these patterns are already covered by
Theorem~\ref{main:classification} and
Proposition~\ref{unified:stabilizer}.
There is an additional simplification of the numerical description:
the largest bounded complementary squares have side length $1/n$.
Consequently, an isometry between two such carpets with parameters
$n$ and $m$ forces $n=m$. Their equivalence is therefore determined by
the $D_4$ orbit of $H$ on this common grid. This simplification uses the
unit-cell hypothesis; for polyomino holes one must in general group
generations as in Theorem~\ref{main:classification}.

For $n=7$, consider the three deletion patterns
\[
\begin{gathered}
 H_1=\{(1,2)\},\qquad H_2=\{(1,2),(5,4)\},\\
 H_4=\{(1,2),(4,1),(5,4),(2,5)\}.
\end{gathered}
\]
Their stabilizers are respectively the trivial group, the group
generated by the half-turn, and the group generated by the quarter-turn.
Indeed, relative to the center of $Q$, the first hole has center vector
proportional to $(-2,-1)$; the other patterns are its orbits under the
half-turn and quarter-turn. Checking the four reflections shows that
none preserves these sets. Proposition~\ref{unified:stabilizer} gives
the same descriptions for their quasisymmetric groups.

\begin{figure}[htbp]
\centering
\begin{tikzpicture}[x=0.34cm,y=0.34cm]
\foreach \offset/\lab in {0/1,10/2,20/4}{
 \begin{scope}[shift={(\offset,0)}]
 \fill[gray!30] (0,0) rectangle (7,7);
 \draw[step=1,gray!60,very thin] (0,0) grid (7,7);
 \draw (0,0) rectangle (7,7);
 \filldraw[fill=white,draw=black,thin] (1,2) rectangle (2,3);
 \ifnum\lab>1
 \filldraw[fill=white,draw=black,thin] (5,4) rectangle (6,5);
 \fi
 \ifnum\lab>2
 \filldraw[fill=white,draw=black,thin] (4,1) rectangle (5,2);
 \filldraw[fill=white,draw=black,thin] (2,5) rectangle (3,6);
 \fi
 \node[below] at (3.5,-0.3) {$H_{\lab}$};
 \end{scope}}
\end{tikzpicture}
\caption{Separated unit-cell deletion patterns with $n=7$.
White cells are holes. Their stabilizers are, from left to right,
the trivial group, the half-turn group, and the quarter-turn group;
none contains a reflection.}
\label{english:patternfigure}
\end{figure}
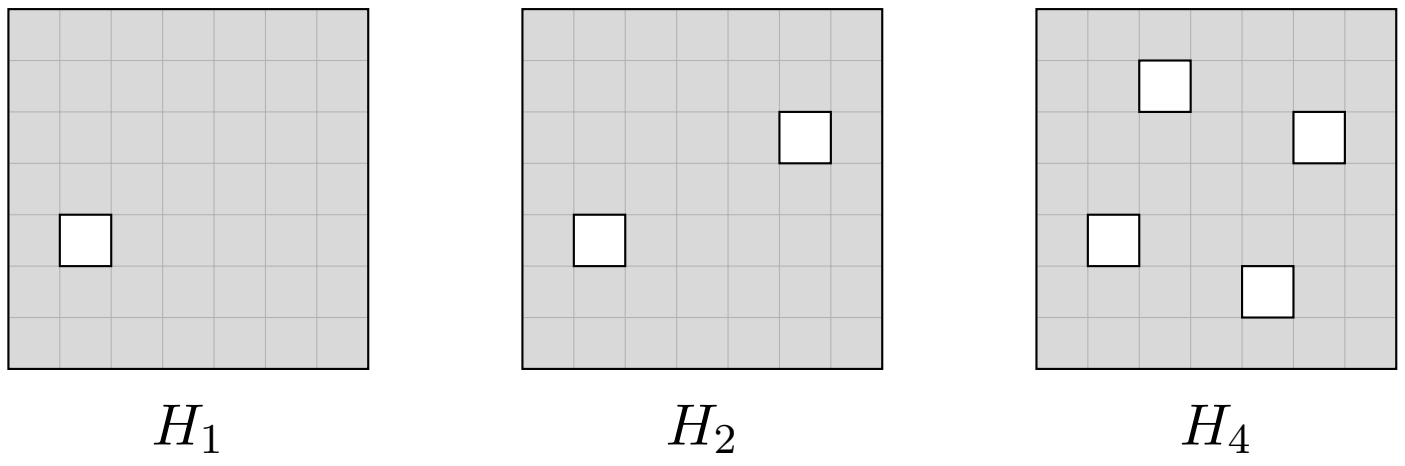

\subsection{An L-shaped hole}

\begin{example}\label{unified:Lexample}
Let $n\ge4$ and $1\le a\le n-3$. Delete exactly the three interior
digits
\[
 (a,a),\qquad(a+1,a),\qquad(a,a+1),
\]
and repeat this pattern with the same orientation in every retained
cell. The three cells merge into one L-shaped actual hole. Then
\[
 \QS(K)=\BiLip(K)=\Iso(K)
 =\{\id|_K,\rho|_K\},\qquad \rho(x,y)=(y,x).
\]
\end{example}

\begin{proof}
The bounds on $a$ place the three deleted cells strictly inside $Q$.
Their merged domain is a Jordan polyomino, so~\eqref{main:P} holds.
The reflection $\rho$ fixes one deleted cell and exchanges the other
two. Conversely, any isometry preserving the hole must fix its unique
reentrant vertex $v=((a+1)/n,(a+1)/n)$ and preserve the missing northeast
quadrant at this vertex. Among the eight linear parts of square
symmetries, only the identity and the coordinate interchange preserve
this quadrant. The corresponding symmetries are precisely $\id$ and
$\rho$, and both fix $v$. Apply
Proposition~\ref{unified:stabilizer}.
\end{proof}

\begin{figure}[htbp]
\centering
\begingroup
\def\CarpetLDigits{0/0,1/0,2/0,3/0,0/1,3/1,0/2,2/2,3/2,0/3,1/3,2/3,3/3}
\newcommand{\CarpetLHole}{%
  \path[fill=white,draw=black,line width=0.22pt]
    (0.25,0.25) -- (0.75,0.25) -- (0.75,0.50)
    -- (0.50,0.50) -- (0.50,0.75) -- (0.25,0.75) -- cycle;}
\begin{tikzpicture}[x=4.2cm,y=4.2cm]
  \begin{scope}
    \fill[gray!28] (0,0) rectangle (1,1);
    \CarpetLHole
    \draw[line width=0.5pt] (0,0) rectangle (1,1);
    \draw[blue!75!black,line width=0.5pt,dash pattern=on 3pt off 2pt]
      (0,0) -- (1,1);
    \node[font=\small,below] at (0.5,-0.045)
      {(a) First generation $Q_1$};
  \end{scope}
  \begin{scope}[shift={(1.45,0)}]
    \fill[gray!28] (0,0) rectangle (1,1);
    \CarpetLHole
    \foreach \i/\j in \CarpetLDigits {
      \begin{scope}[shift={(\i/4,\j/4)},scale=0.25]
        \CarpetLHole
        \foreach \k/\l in \CarpetLDigits {
          \begin{scope}[shift={(\k/4,\l/4)},scale=0.25]
            \CarpetLHole
          \end{scope}
        }
      \end{scope}
    }
    \draw[line width=0.5pt] (0,0) rectangle (1,1);
    \draw[blue!75!black,line width=0.5pt,dash pattern=on 3pt off 2pt]
      (0,0) -- (1,1);
    \node[font=\small,below] at (0.5,-0.045)
      {(b) Third approximation $Q_3$};
  \end{scope}
\end{tikzpicture}
\endgroup
\caption{The L-shaped deletion pattern with $n=4$ and $a=1$.
The deleted first-level digits are $(1,1)$, $(2,1)$, and $(1,2)$;
their union forms one actual hole. The right panel is the third finite
approximation, with $1+13+13^2=183$ L-shaped holes from the first
three generations. Gray denotes the retained region, and the blue
dashed line marks the reflection axis $y=x$ in each panel.}
\label{unified:Lfigure}
\end{figure}
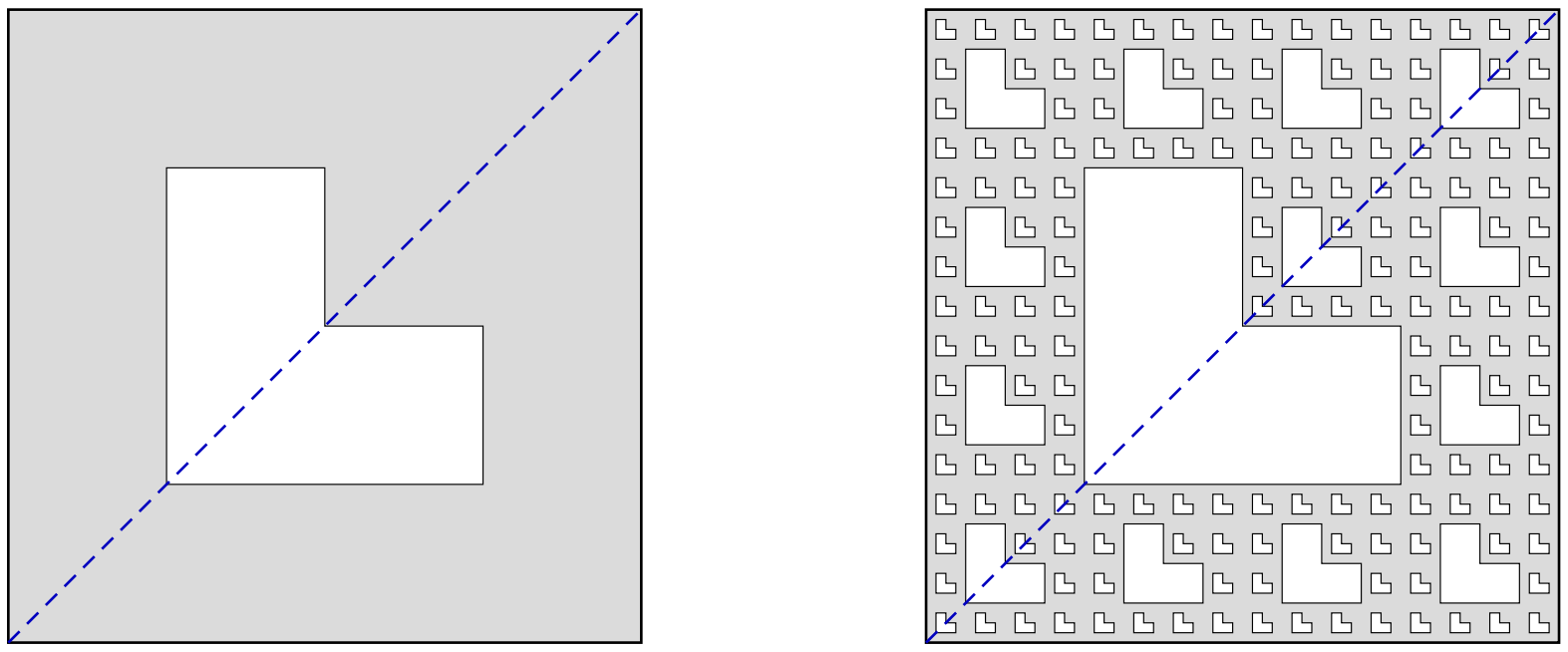

\subsection{Changing the four marked points on the outer square}
\label{scope:remarking}

\begin{example}\label{scope:markedexample}
There is an area-zero Sierpi\'nski carpet $T$ with outer square $Q$,
axis-parallel square holes, and uniformly relatively separated
peripheral circles, for which
\[
       \QS(T)\cong D_4,\qquad \Iso(T)\cong C_4.
\]
A quasisymmetric involution of $T$ takes a corner to the interior
of a side. The carpet $T$ is quasisymmetrically equivalent to the
standard middle-third carpet $S_3$, but the two are not isometric.
\end{example}

\begin{proof}
Let $r(x,y)=(1-y,x)$ be the counterclockwise quarter-turn of $Q$.
Starting at $p_0=(1/3,0)$, mark $p_j=r^j(p_0)$, $0\le j\le3$,
on the outer boundary of $S_3$. The four points, in positive order, are
\[
 (1/3,0),\quad (1,1/3),\quad (2/3,1),\quad (0,2/3).
\]
The carpet $S_3$ has area zero since the area of its $k$th
approximation is $(8/9)^k$. Its peripheral circles are squares
and have relative separation at least $1/\sqrt{2}$. Indeed, the
smaller of two peripheral squares, of side $s$, is at distance $s$
from the boundary of its parent cell. The other peripheral circle
lies outside the interior of that cell. This argument also applies
when the other circle is $\partial Q$.

Use the four successive boundary arcs as the sides in Ntalampekos's
square uniformization. In the numbering of the preprint
\cite[Chap.~3, Thm.~3.2, Cor.~3.49 and Cor.~3.60]{Nta20}, the
choice of sides is made in \S4.2 and their separation is used in
\S11, (3.56). Here $\min_{i\ne j}|p_i-p_j|=\sqrt{5}/3$.
We obtain an orientation-preserving quasisymmetry
$\phi:S_3\to T$ onto an area-zero carpet in $[0,1]\times[0,H]$,
with axis-parallel square holes. Its outer corners
$b_j=\phi(p_j)$ start at $b_0=(0,0)$ and proceed counterclockwise,
with $b_1=(1,0)$. Each of the four arcs maps onto the corresponding
side. The peripheral circles of $T$ are uniformly relatively
separated because $\phi$ is quasisymmetric.

We first determine the height $H$. The conjugate
$h=\phi r\phi^{-1}$ sends $b_j$ to $b_{j+1}$, with indices taken
modulo four. Write $J(x,y)=(y,1-x)$ and set
\[
       G:H^{-1}T\longrightarrow J(T),\qquad G(z)=J(h(Hz)).
\]
The two outer rectangles have height one and widths $1/H$ and $H$,
respectively. The map $G$ preserves orientation, sends corners to
corners, and fixes the origin, since $J(b_1)=(0,0)$.
The four-corner rigidity theorem of Bonk and Merenkov
\cite[Thm.~1.4]{BM2013} gives $1/H=H$ and $G=\id$ on its domain.
Thus $H=1$ and $h=J^{-1}|_T=r|_T$.

Now take the reflection $\sigma(x,y)=(1-x,y)$ of $S_3$ and set
$f=\phi\sigma\phi^{-1}$. The point $\sigma(p_0)=(2/3,0)$ lies on
the open counterclockwise boundary arc from $p_0$ to $p_1$. Hence
$f(b_0)$ belongs to the open bottom side of $Q$
(Figure~\ref{scope:markfigure}). Every point of that side is at
distance strictly less than $\sqrt{2}$ from every point of $Q$.
In particular,
\[
              |f(b_0)-f(b_2)|<\sqrt{2}=|b_0-b_2|.
\]

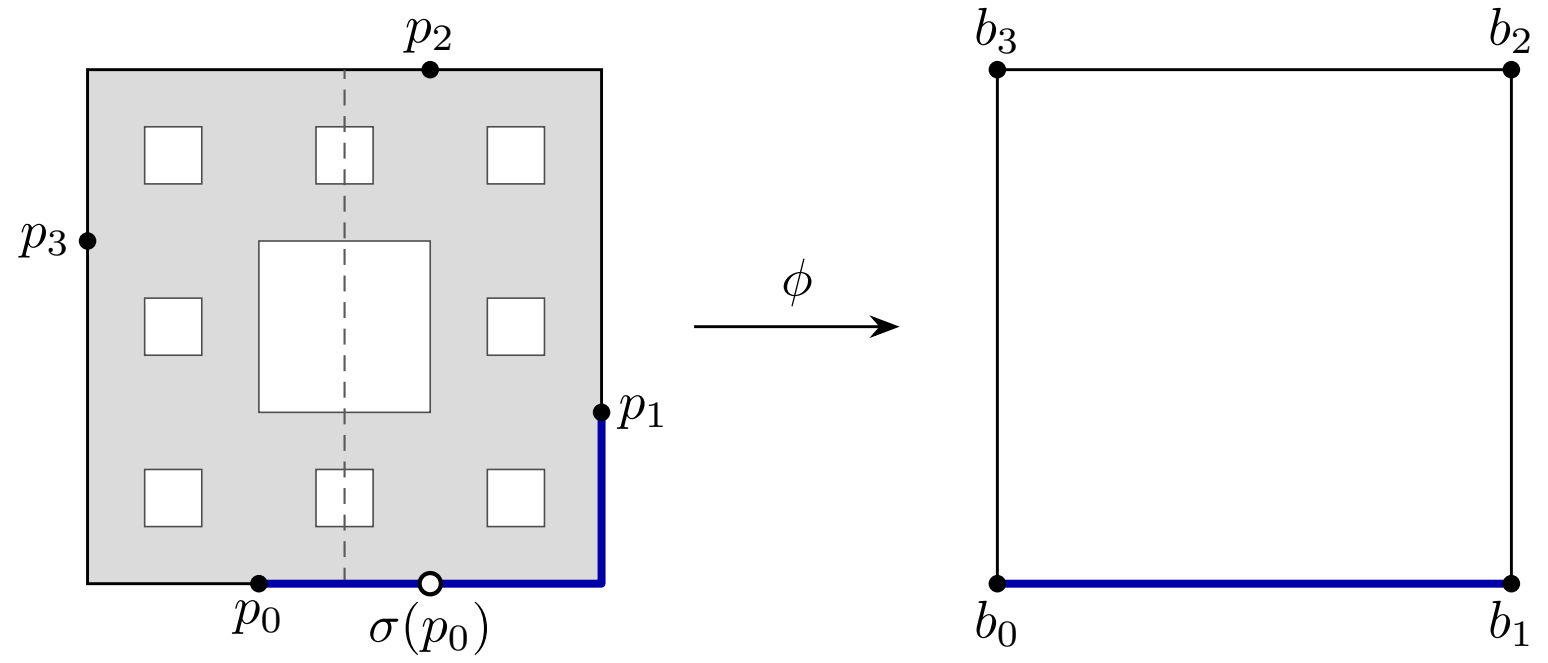
\begin{figure}[htbp]
\centering
\begin{tikzpicture}[x=3.4cm,y=3.4cm,font=\small,
    marked arc/.style={draw=blue!65!black,line width=1.5pt,line join=round},
    mark label/.style={inner sep=2pt}]
  \fill[gray!28] (0,0) rectangle (1,1);
  \path[draw=black!70,fill=white,line width=.3pt]
    (1/3,1/3) rectangle (2/3,2/3);
  \foreach \i/\j in {0/0,1/0,2/0,0/1,2/1,0/2,1/2,2/2}{
    \path[draw=black!70,fill=white,line width=.3pt]
      ({(3*\i+1)/9},{(3*\j+1)/9})
      rectangle ({(3*\i+2)/9},{(3*\j+2)/9});
  }
  \draw[line width=.55pt] (0,0) rectangle (1,1);
  \draw[densely dashed,black!65,line width=.4pt] (.5,0) -- (.5,1);
  \draw[marked arc] (1/3,0) -- (1,0) -- (1,1/3);
  \foreach \x/\y/\i/\pos in
      {0.333333/0/0/below,1/0.333333/1/right,
       0.666667/1/2/above,0/0.666667/3/left}{
    \fill (\x,\y) circle[radius=1.65pt]
      node[mark label,\pos=1pt] {$p_\i$};
  }
  \draw[fill=white,line width=.8pt] (2/3,0) circle[radius=2pt]
    node[mark label,below=1pt] {$\sigma(p_0)$};
  \draw[-{Stealth[length=2mm]},line width=.55pt]
    (1.18,.5) -- (1.58,.5) node[midway,above] {$\phi$};
  \begin{scope}[shift={(1.77,0)}]
    \draw[line width=.55pt] (0,0) rectangle (1,1);
    \draw[marked arc] (0,0) -- (1,0);
    \foreach \x/\y/\i/\pos in
        {0/0/0/below,1/0/1/below,1/1/2/above,0/1/3/above}{
      \fill (\x,\y) circle[radius=1.65pt]
        node[mark label,\pos=1pt] {$b_\i$};
    }
  \end{scope}
\end{tikzpicture}
\caption{Changing the four marked points on the outer boundary.
On the second approximation to $S_3$ (left), reflection in the dashed
line sends $p_0$ to the open circle. The highlighted boundary arc from
$p_0$ to $p_1$ maps to the bottom side of $T$ (right, outer boundary
only). Thus $\phi(\sigma(p_0))$ lies strictly between $b_0$ and $b_1$.}
\label{scope:markfigure}
\end{figure}

Bonk and Merenkov \cite[Thm.~1.1]{BM2013} proved that
$\QS(S_3)=D_4|_{S_3}$. Conjugation by $\phi$ gives all eight elements
of $\QS(T)$. Four are the Euclidean rotations $h^j$, $0\le j<4$.
Each remaining element $h^jf$ shortens the same diagonal pair as $f$.
Consequently $\Iso(T)=\langle h\rangle$. An isometry from $S_3$
onto $T$ would conjugate $\Iso(S_3)$ to $\Iso(T)$, which is
impossible since their orders are eight and four.
\end{proof}

There is also an exact description of how $f$ acts on holes at
every generation. Let $\mathcal H_k$ be the family of holes created
at stage $k$ of $S_3$. To each $V\in\mathcal H_k$ associate the
unique open square $U_V$ in the complement of $T$ whose boundary
is $\phi(\partial V)$. The uniformization gives every inner hole
of $T$ exactly once in this way. Thus
\[
 T=Q\setminus\bigcup_{k\ge1}\ \bigcup_{V\in\mathcal H_k}U_V,
 \qquad f(\partial U_V)=\partial U_{\sigma(V)}.
\]
The second identity follows from $f\phi=\phi\sigma$. Since
$\sigma$ preserves the ternary grid and its deletion rule,
$\sigma(\mathcal H_k)=\mathcal H_k$ for every $k$. This indexes the
target holes in families of $8^{k-1}$. The map $\phi$ determines
their sizes and positions.

\section{Scope of the results}

Example~\ref{scope:markedexample} separates the geometry of the
peripheral squares from the repetition of a fixed grid pattern.
Its target $T$ cannot have a presentation satisfying
Theorem~\ref{main:classification}, since it admits the
corner-moving involution constructed above.

Repetition of the pattern supplies the same tail carpet in each cell
and makes the periodic extension forward invariant under integer
expansion. Grouping over a period preserves these properties.
A nonperiodic initial block or an arbitrary sequence of patterns
requires a separate argument. For example, for central holes of
varying relative side lengths $a_k$, Mackay, Tyson and
Wildrick~\cite{MTW13} proved that their carpets with the natural
measure support a $1$-Poincar\'e inequality precisely when
$\sum_k a_k<\infty$. Positive area lies outside the measure-zero
rigidity argument used here. For touching hole closures or local
cut points, the common-side propagation also requires a different
analysis.

\end{document}